\documentclass{article}

\usepackage[utf8]{inputenc}
\usepackage{lmodern}
\usepackage[margin=2.5em,font=small]{caption}
\usepackage{macros}
\usepackage{math}
\usepackage[section]{theorems}
\usepackage{animals}
\usepackage{enumitem}
\setlist{itemsep=0em}

\usetikzlibrary{decorations.pathreplacing}
\tikzstyle{arrow}=[-stealth,very thick]
\tikzstyle{brace}=[very thick, decorate, decoration=brace]

\title{Average site perimeter of directed animals on the two-dimensional
lattices}
\author{Axel Bacher}
\date\today

\begin{document}
\maketitle

\begin{abstract}
We introduce new combinatorial (bijective) methods that enable us to compute
the average value of three parameters of directed animals of a given area,
including the site perimeter. Our results cover directed animals of any
one-line source on the square lattice and its bounded variants, and we give
counterparts for most of them in the triangular lattices. We thus prove
conjectures by Conway and Le~Borgne. The techniques used are based on
Viennot's correspondence between directed animals and heaps of pieces (or
elements of a partially commutative monoid).
\end{abstract}

\section{Introduction}\label{intro}

Let $\Gamma$ be an oriented graph and $S$ a nonempty finite set of vertices of
$\Gamma$. A \emph{directed animal} of source $S$ on $\Gamma$ is a finite set
of vertices $A$ that contains $S$ and such that for every vertex $v$ of $A$,
there exists a vertex $s$ of $S$ and a path from $s$ to $v$ going only through
vertices of $A$. The vertices of a directed animal $A$ are called
\emph{sites}. The \emph{area} of $A$, denoted by $\abs A$, is the number of
sites of $A$.
 
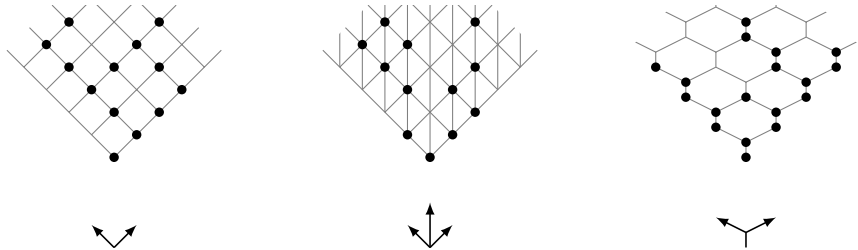
\begin{figure}[ht]
\begin{center}
\begin{tikzpicture}[scale=.3]
\clip (-19,-4.25) rectangle (19,6.75);
\begin{scope}[xshift=-14cm,animal]
\draw[help lines] (0,0) grid (4.75,4.75);
%\draw[-latex] (0,0) -- (5.25,0) node[anchor=base west] {$i$};
%\draw[-latex] (0,0) -- (0,5.25) node[anchor=base east] {$j$};
\foreach \p in {(0,0),(1,0),(2,0),(3,0),(1,1),(3,1),(1,2),(2,2),(3,2),(4,2),
(1,3),(1,4),(2,4)}
\node at \p [site] {};
\draw[marked,-latex] (-2,-2) -- +(0,1);
\draw[marked,-latex] (-2,-2) -- +(1,0);
\end{scope}
\begin{scope}[xshift=0cm]
\begin{scope}[animal] 
%\clip (-.5,-.5) rectangle (4.5,3.5); 
\draw [help lines] (0,0) grid (4.75,4.75); 
%\draw[-latex] (0,0) -- (5.25,0) node[anchor=base west] {$i$};
%\draw[-latex] (0,0) -- (0,5.25) node[anchor=base east] {$j$};
\draw[marked,-latex] (-2,-2) -- +(0,1);
\draw[marked,-latex] (-2,-2) -- +(1,0);
\draw[marked,-latex] (-2,-2) -- +(1,1);
\end{scope}
\begin{scope}[animal]
\clip (-.5,-.5) rectangle (4.75,4.75);
\draw[help lines] (0,4) -- ++(4,4) (0,3) -- ++(4,4) (0,2) -- ++(4,4)
(0,1) -- ++(3.5,3.5) (0,0) -- ++(4,4) (1,0) -- ++(4,4)
(2,0) -- ++(4,4) (3,0) -- ++(4,4) (4,0) -- ++(4,4);
\foreach \p in {(0,0),(1,0),(2,0),(0,1),(2,1),(3,1),(1,2),(2,4),(1,4),(4,2),(1,3),
(2,3)}
\node at \p [site] {};
\end{scope}
\end{scope}
\begin{scope}[xshift=14cm,scale=1.333,yshift=.5cm]
\begin{scope}
\clip[animal] (-.5,-.5) rectangle ++(4,4);
\foreach \p in
{(0,0),(-1,1),(1,1),(-2,2),(0,2),(2,2),(-3,3),(-1,3),(1,3),(3,3),(-2,4),(0,4),(2,4)}
\draw[help lines] \p -- +(0,-.5) \p -- +(1,.5) \p -- +(-1,.5);
\draw[help lines] (-1,4.5) -- ++(0,.5) (1,4.5) -- ++(0,.5);
%(-1,.5) -- (0,0) -- (1,.5);
\end{scope}
\foreach \p in {(0,-.5),(0,0),(-1,.5),(1,.5),(-1,1),(1,1),(-2,1.5),(0,1.5),(2,1.5),
(-2,2),(2,2),(-3,2.5),(1,2.5),(3,2.5),(1,3),(3,3),(0,3.5),
(0,4)}
\node at \p [site] {};
\draw[marked,-latex] (0,-3.5) -- ++(0,.5) -- +(1,.5);
\draw[marked,-latex] (0,-3) -- +(-1,.5);
\end{scope}
\end{tikzpicture}
\end{center}
\caption{Single-source directed animals on the square, triangular and
honeycomb lattices. All edges point upwards.}\label{animals}
\end{figure}

On Figure~\ref{animals} are depicted single-source directed animals on the
three two-dimensional regular lattices: the square lattice, the triangular
lattice, and the honeycomb lattice.

Single-source directed animals constitute a subclass of \emph{animals} (an
animal on a non-oriented graph $\Gamma$ is simply a finite connected set of
vertices of $\Gamma$).  While the enumeration of animals on any lattice is an
open problem despite extensive research for decades, directed animals are
fairly easier to enumerate. As we will not deal with general animals in this
paper, we will abusively use the term \emph{animal} instead of directed
animal.

Single-source directed animals on the square and triangular lattices have been
enumerated \cite{dhar,gouyou,betrema}. Specifically, let $a(n)$ and $\bar a(n)$
be the number of animals of source $\{(0,0)\}$ and area $n$ on the square and
triangular lattice, respectively. The generating functions of these numbers
are:
\begin{align}\label{introsquare}
\sum_{n\geq1}a(n)t^n&=\fullsquare\text;\\\label{introtriang}
\sum_{n\geq1}\bar a(n)t^n&=\fulltriang.
\end{align}

Even then, much remains unclear. The enumeration of directed animals on the
honeycomb lattice is an open problem, and according to \cite{guttmann}, the
generating function is probably not D-finite; on the square and triangular
lattices, comparatively very little is known when one tries to take into
account parameters other than area.

\bigskip

Today, two enumeration methods account for almost every known result on
directed animals. One of them is the \emph{gas model} technique, originally
used by Dhar \cite{dhar}. This technique was further developed by
Bousquet-Mélou \cite{bousquet}; see also \cite{marckert,albenque} for more
recent work.

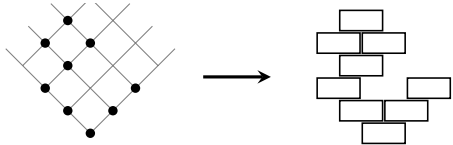
\begin{figure}[ht]
\centering
\begin{tikzpicture}[scale=.3]
\clip (-10.25,-.5) rectangle (10.25,5.75);
\draw[arrow] (-1.5,2.5) -- (1.5,2.5);
\begin{scope}[xshift=-6.5cm,animal]
\draw[help lines] (0,0) grid (3.75,3.75); 
\foreach \p in {(0,0),(1,0),(2,0),(0,1),(0,2),(1,2),(2,2),(1,3),(2,3)} 
\node at \p [site] {};
\end{scope}
\begin{scope}[reset cm,dominoes,xshift=6.5cm,animal]
\foreach \p in {(0,0),(1,0),(2,0),(0,1),(0,2),(1,2),(2,2),(1,3),(2,3)} 
\node at \p [domino] {};
\end{scope}
\end{tikzpicture}
\caption{A directed animal on the square lattice can be turned into a heap by
replacing each site by a $2{\times}1$ domino.}\label{animal-heap}
\end{figure}

The method used in this paper is the second one, based on a correspondence,
due to Viennot \cite{viennot}, between animals and other objects called
\emph{heaps of dominoes}. The basic idea is to replace each site of an animal
by a $2{\times}1$ domino, so that each domino either lies on the ground or
sits on one or two other dominoes (Figure~\ref{animal-heap}).

As we will see later, this method works for the triangular lattice as well.
However, no simple model of heaps of dominoes has been found to correspond to
animals on the honeycomb lattice. This may explain the lack of knowledge on
the subject.

\bigskip

\begin{figure}[ht]
\centering
\begin{tikzpicture}[scale=.3]
\draw [-latex] (0,-2) -- +(6,6) node [anchor=west] {$i$};
\draw [-latex] (0,-2) -- +(-6,6) node [anchor=east] {$j$};
\begin{scope}
\clip (-6,-2.5) rectangle (6,6.75);
\begin{scope}[animal]
\draw[help lines] (-1,-1) grid (4.75,4.75);
\foreach \p in {(0,0),(1,0),(2,0),(1,1),(2,1),(3,1),(3,2),(1,2),(3,3),
(1,3),(1,4),(2,4),(4,1),(-1,-1),(-1,0),(-1,1)}
\node at \p [site] {};
\draw[marked,rounded corners=.12cm] (3.4,3) -- ++(-.4,-.4) -- ++(-1.4,1.4)
-- ++(.4,.4) -- cycle;
\draw[marked,rounded corners=.12cm] (1.717,.717) |- ++(.566,-1)
|- ++(-1.566,1.566) -- ++(0,-.566) [sharp corners] -- cycle;
\draw[marked] (-1,2) circle (.283);
\end{scope}
\end{scope}
\begin{scope}[xshift=10cm,animal]
\draw[marked,rounded corners=.12cm] (3.4,2) -- ++(-.4,-.4) -- ++(-1.4,1.4) -- ++(.4,.4) -- cycle;
\draw[marked,rounded corners=.12cm] (.717,.717) |- ++(.566,-1)
|- ++(-1.566,1.566) -- ++(0,-.566) [sharp corners] -- cycle;
\draw[marked] (-.5,-.5) circle (.283);
\node at (3.5,1.5) [anchor=west] {adjacent sites};
\node at (2,0) [anchor=west] {loop};
\node at (.5,-1.5) [anchor=west] {neighbour};
\end{scope}
\end{tikzpicture}
\caption{A directed animal on the square lattice with two adjacent sites, a
loop, and a neighbour marked.}\label{parameters}
\end{figure}
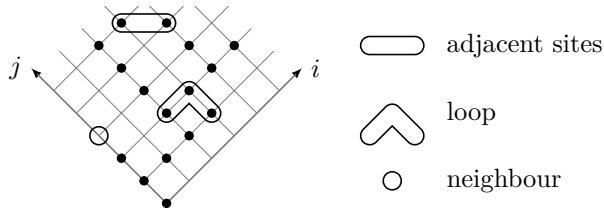

The purpose of this paper is to study three other parameters of directed
animals, introduced below, and illustrated in Figure~\ref{parameters}:
\begin{itemize}
\item two sites of an animal on the square or triangular lattice are
\emph{adjacent} if they are of the form $(i+1,j)$ and $(i,j+1)$. We denote by
$j(A)$ the number of pairs of adjacent sites of $A$.
\item a \emph{loop} consists of two adjacent sites $(i+1,j)$ and $(i,j+1)$,
along with a third site at $(i+1,j+1)$. We denote by $\ell(A)$ the number of
loops of $A$.
\item a \emph{neighbour} of an animal $A$ of source $S$ is a vertex $v$ not in
$A$, such that $A\cup\{v\}$ is still a directed animal of source $S$. The
number of neighbours of $A$ is called the \emph{site perimeter} of $A$ and is
denoted by $p(A)$.
\end{itemize}

Taking, for instance, the site perimeter, we may consider the bivariate
generating function counting single-source animals according to both area and
perimeter on the square lattice:
\[A^p(t,u)=\sum_At^\abs Au^{p(A)}.\]
This generating function is not known, and is believed not to be D-finite
\cite{jensen}. Instead, we will consider the generating function giving the
total number of neighbours in the animals of a fixed area:
\[\sum_Ap(A)t^\abs A=\diff{A^p}{u}(t,1).\]
By dividing the total site perimeter of the animals of area $n$ by the number
of these animals, one gets the \emph{average} site perimeter in animals of a
fixed area. Alternatively, this generating function may be seen as counting
single-source directed animals with a marked neighbour.

This function, and the ones that similarly give the average number of adjacent
sites and loops, turns out to be easier to derive. Specifically, the value of
the generating function counting the total number of loops of single-source
animals on the square lattice was obtained by Bousquet-Mélou using gas model
methods \cite{bousquet}:
\begin{equation}\label{introloops}
\sum_A\ell(A)t^\abs A=\loops.
\end{equation}
As for the total site perimeter on the square lattice, it was the object of a
conjecture by Conway in 1996 \cite{conway}:
\begin{equation}\label{introperimeter}
\sum_Ap(A)t^\abs A=\conway.
\end{equation}
Le Borgne \cite{leborgne} also conjectured the value of similar generating
functions counting the site perimeter of animals on square and triangular
lattices of bounded width.

In Section~\ref{section:results}, we prove these conjectures and give a new
proof of \eqref{introloops} using combinatorial methods; moreover, we show
that the total number of adjacent sites is given by:
\begin{equation}\label{introadjacent}
\sum_Aj(A)t^\abs A=\adjacent.
\end{equation}
Actually, our results are more general than that: the same methods can be used
on different kinds of lattices, obtained by adding one or two vertical walls
(the \emph{half-lattice}, \emph{cylindrical lattices} and \emph{rectangular
lattices}, defined in Section~\ref{section:dominoes}), and on animals with any
fixed source.
  
Knowing, say, the total site perimeter of single-source animals of area~$n$ on
the square lattice, we get their average perimeter by dividing by the number
$a(n)$ of these animals:
\begin{equation*}
p(n)=\frac1{a(n)}\sum_{\abs A=n}p(A).
\end{equation*}
This quantity may thus be computed using \eqref{introsquare} and
\eqref{introperimeter}. The numbers of adjacent sites and loops are handled
similarly. From these generating functions, singularity analysis
\cite{flajolet} yields estimates on these quantities as $n$ tends to infinity:
\begin{align*}
j(n)&\sim\frac n4;&\ell(n)&\sim\frac n9;&p(n)\sim\frac{3n}4.
\end{align*}

\bigskip

The paper is organized as follows. In Section~2, we introduce in detail the
notion of heaps of pieces and give several lemmas useful for animal
enumeration. In Section~3, we enumerate directed animals of any source on
several kinds of square and triangular lattices, according to area alone.
In Section~4, we give a general method to derive the generating functions
giving the average number of adjacent sites, number of loops, and site
perimeter of directed animals on the square lattice, as well as counterparts 
of most of these results on the triangular lattice. We derive asymptotic
results in Section~5. Finally, we illustrate our formul{\ae} with a few
examples in Section~6.

\section{Heaps of pieces}\label{section:heaps}

The notion of heaps of pieces is due to Viennot, and this topic is covered in
detail in \cite{viennot}. We repeat the definitions for convenience, and make
a few minor additions which we use later.

\subsection{Basics}

Intuitively, a heap is a finite set of pieces. It is built by dropping
successively the pieces at certain positions, chosen from a given set. When
the positions of two pieces overlap, the second piece falls on the first, like
in Figure~\ref{animal-heap}. A formal definition is given below.

\begin{definition}\label{definition:heap}
Let $Q$ be a set and $\cC$ a reflexive symmetric relation on $Q$. A
\emph{heap} of the model $(Q,\cC)$ is a finite subset $H$ of
$Q\times\bN$ satisfying:
\begin{enumerate}
\item if $(q,i)$ and $(q',i)$ with $q\neq q'$ are in $H$, then $(q,q')$ is not
in $\cC$;
\item if $(q,i)$ is in $H$ and $i>0$, then there exists $(q',i-1)$ in $H$ such
that $(q,q')$ is in $\cC$.
\end{enumerate}
\end{definition}

The relation $\cC$ is called the \emph{concurrency relation}, and two
positions $q$ and $q'$ are \emph{concurrent} if $(q,q')$ is in $\cC$ (in the
above intuitive definition, this means that they overlap). The elements of a
heap are called \emph{pieces}. If $(q,i)$ is a piece of a heap, $q$ is called
its \emph{position} and $i$ its \emph{height}.

The pieces of a heap are naturally equipped with a poset structure: define the
relation $\prec$ such that  $(q,i)\prec(q',i')$ whenever $q$ and $q'$ are
concurrent and $i<i'$. Let $\leq$ be the reflexive transitive closure of
$\prec$. We say that a piece $x$ is \emph{below} a piece $y$, or $y$ is
\emph{above} $x$, if $x\leq y$. The pieces $x$ and $y$ are \emph{independent}
if neither is above the other.

The partial order $\leq$ may be viewed more intuitively: let $H$ be a heap and
$x$ a piece of $H$. If one takes the piece $x$ and push it upwards, it pushes
along some pieces in the way. If one pushes it high enough, the moved pieces
are exactly the pieces above $x$.

The pieces of a heap $H$ that are minimal for the order $\leq$ are called
\emph{minimal pieces}; they are exactly the pieces of height~0. The set of
their positions is called the \emph{base} of $H$, and denoted by $b(H)$.
Likewise, the pieces that are maximal for $\leq$ are called \emph{maximal
pieces}; we denote by $m(H)$ the set of their positions.

We now define generating functions counting the heaps of a model; we denote by
$\abs H$ the number of pieces of the heap $H$.

\begin{definition}\label{definition:gf1}
Let $(Q,\cC)$ be a model of heaps and $S$ a finite subset of $Q$. We denote by
$\cH_S(t)$ and $\cH_{[S]}(t)$ the generating functions (provided they exist)
of heaps respectively of base $S$ and with base included in $S$:
\begin{align*}
\cH_S(t)&=\sum_{b(H)=S}t^{\abs H}\text;&
\cH_{[S]}(t)&=\sum_{b(H)\subseteq S}t^{\abs H}\text.
\end{align*}
\end{definition}

These two generating functions are obviously linked by:
\begin{align*}
\cH_{[S]}(t)&=\sum_{T\subseteq S}\cH_T(t)\text;
\intertext{Conversely, the inclusion-exclusion principle yields:}
\cH_S(t)&=\sum_{T\subseteq S}(-1)^{\abs {S\smallsetminus T}}\cH_{[T]}(t)\text.
\end{align*}

Let us now assume that the set of positions $Q$ of the model is finite. In
this case, the generating functions above may be computed using a result due
to Viennot \cite{viennot}, which we present below.

A heap is called \emph{trivial} if all its pieces have height~0. This means
that a trivial heap may be identified with the set of the positions of its
pieces, which must be pairwise nonconcurrent. As $Q$ is a finite set, there is
only a finite number of trivial heaps. We denote by $\cT_{[S]}(t)$ the
alternating generating function of trivial heaps included in $S$:
\[\cT_{[S]}(t)=
\sum_{\substack{T\subseteq S\text,\\T\text{ trivial}}}(-t)^{\abs T}\text.\]
Since $Q$ is finite, this generating function is actually a polynomial,
and is usually relatively easy to compute.

\begin{lemma}[Inversion Lemma]\label{lemma:inversion}
Let $(Q,\cC)$ be a finite heap model and $S$ a subset of $Q$. The generating
function of heaps of base included in $S$ is:
\[\cH_{[S]}(t)=\frac{\cT_{[Q\smallsetminus S]}(t)}{\cT_{[Q]}(t)}.\]
\end{lemma}

Thanks to this lemma, we see that in a finite model, the generating function
$\cH_{[S]}(t)$ is a quotient of two polynomials, and is therefore rational.

\subsection{Strict and inflated heaps}

In this section, we define families of heaps which we use in the
correspondence with directed animals. Let $H$ be a heap, and let $(q,i)$ and
$(q',i')$ be two pieces. We say that $(q',i')$ \emph{sits} on $(q,i)$ if $q$
and $q'$ are concurrent and $i'=i+1$. Thus, Condition~2 of
Definition~\ref{definition:heap} states that any non-minimal piece must sit on
another piece.

The objects obtained by relaxing this condition, keeping only Condition~1, are
called \emph{pre-heaps}.

\begin{definition}
A heap or pre-heap is \emph{strict} if it has no piece sitting on another at
the same position, i.e.\ no two pieces $(q,i)$ and $(q,i+1)$.
\end{definition}

\begin{definition}
An \emph{inflated heap} is a strict pre-heap $H$, such that for every piece
$(q,i)$ satisfying $i>0$, at least one of the following pieces is in $H$:
\begin{enumerate}[label=(\alph*)]
\item a piece $(q',i-1)$, such that $q$ and $q'$ are concurrent (and $q\neq
q'$);
\item the piece $(q,i-2)$.
\end{enumerate}
\end{definition}

Examples are found in Figure~\ref{figure:models}.

Let $H$ be a heap. A \emph{stack} of $H$ is a maximal set of pieces all at the
same position, with consecutive heights. Thus, a heap is strict if all its
stacks have only one piece.

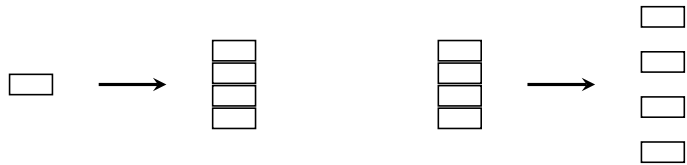
\begin{figure}[ht]
\centering
\begin{tikzpicture}[dominoes]
\node at (0,-1.5) [domino] {};
\node at (0,-.5) [domino] {};
\node at (0,.5) [domino] {};
\node at (0,1.5) [domino] {};
\node at (9,-3) [domino] {};
\node at (9,-1) [domino] {};
\node at (9,1) [domino] {};
\node at (9,3) [domino] {};
%\draw[help lines] (1,-2) -- (4,-3.5) (1,2) -- (4,3.5);
\draw[arrow] (3,0) -- ++(3,0);
\begin{scope}[xshift=-19cm]
\node at (0,0) [domino] {};
\node at (9,-1.5) [domino] {};
\node at (9,-.5) [domino] {};
\node at (9,.5) [domino] {};
\node at (9,1.5) [domino] {};
%\draw[help lines] (1,-.5) -- (4,-2) (1,.5) -- (4,2);
\draw[arrow] (3,0) -- ++(3,0);
\end{scope}
\end{tikzpicture}
\caption{By replacing each piece of a strict heap by a stack of pieces, one
gets a general heap; by inflating each stack, one gets an inflated heap.}
\label{figure:stack}
\end{figure}

Any heap may be built in a unique manner from a strict heap by replacing each
piece by a stack consisting of an arbitrary positive number of pieces; in
turn, any inflated heap may be built from a heap by ``inflating'' each stack,
pushing along all pieces above it (see Figure~\ref{figure:stack}).

These remarks enable us to derive from Lemma~\ref{lemma:inversion} the
generating functions of strict and inflated heaps. First, as inflating a stack
does not change its base or number of pieces, the generating functions of
inflated heaps are the same as those of general heaps.

\begin{notation}\label{definition:gf2}
Let $(Q,\cC)$ be a model of heaps. We denote by $H_S(t)$ and $H_{[S]}(t)$ the
generating functions of strict heaps with base $S$ and base included in $S$,
respectively%
\footnote{%
For the sake of clarity, all generating functions in this paper follow the
same typographical pattern as Definitions \ref{definition:gf1} and
\ref{definition:gf2}: the generating functions of general (or inflated) heaps
are denoted by calligraphic letters, while the ones of strict heaps are
denoted by standard capital letters. Likewise, a subscript $[S]$ is square
brackets always indicates heaps with a base included in $S$, while a subscript
$S$ denotes heaps of base $S$.}%
.
\end{notation}

The construction of Figure~\ref{figure:stack} translates into a link between
the generating functions of strict and general heaps:
\begin{align*}
\cH_S(t)&=H_S\biggl(\frac{t}{1-t}\biggr)\text,
\intertext{or equivalently:}
H_S(t)&=\cH_S\biggl(\frac{t}{1+t}\biggr)\text.
\end{align*}
These links remain valid between generating functions of heaps with base
included in~$S$.

\subsection{Factorized heaps}

We now present a monoid structure, again due to Viennot, on the set of heaps
of a given model. Let $H$ be a heap and $q$ a position. Let $H\cdot q$ be the
heap formed by dropping a piece at position $q$ on top of $H$; more formally,
let $H\cdot q$ be $H\cup\{(q,i)\}$, where $i$ is the largest integer such that
this is a heap.

In this way, a heap may be built one piece at a time. This may be done by
adding the pieces in any order compatible with the partial order $\leq$. This
idea is used to define the product of two heaps.

\begin{definition}\label{definition:product}
Let $H_1$ and $H_2$ be two heaps. The product $H_1\cdot H_2$ is built by
letting all pieces of $H_2$ fall on $H_1$, in any order compatible with
$\leq$.

A \emph{factorized heap} is a heap $H$, with a distinguished factorization
$H=H_1\cdot H_2$. We denote such a heap $(H=H_1\cdot H_2)$ or
$(H_1\cdot H_2)$. A factorized heap $(H=H_1\cdot H_2)$ is \emph{strict} if $H$
is strict; it is \emph{almost strict} if both $H_1$ and $H_2$ are strict.
\end{definition}

\begin{figure}[ht]
\centering
\begin{tikzpicture}[dominoes]
\begin{scope}
\node at (0,0) [marked domino] {};
\node at (-1,1) [marked domino] {};
\node at (-2,2) [marked domino] {};
\node at (-1,3) [marked domino] {};
\node at (-2,4) [marked domino] {};
\node at (1,6) [domino] {};
\node at (0,7) [domino] {};
\node at (2,7) [domino] {};
\node at (-1,8) [domino] {};
\end{scope}
\draw[arrow] (5,3) -- ++(3,0);
\begin{scope}[xshift=13cm]
\node at (0,0) [marked domino] {};
\node at (-1,1) [marked domino] {};
\node at (-2,2) [marked domino] {};
\node at (-1,3) [marked domino] {};
\node at (-2,4) [marked domino] {};
\node at (1,1) [domino] {};
\node at (0,4) [domino] {};
\node at (2,2) [domino] {};
\node at (-1,5) [domino] {};
\end{scope}
\end{tikzpicture}
\caption{The product of two heaps is obtained by dropping the second heap on
top of the first.}
\label{figure:product}
\end{figure}
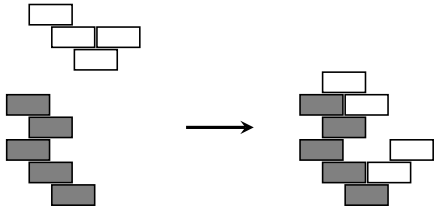

The monoid structure induced by this product is isomorphic to the
\emph{partially commutative monoid} \cite{cartier} on the alphabet $Q$, with
concurrency relation $\cC$. The product is illustrated in
Figure~\ref{figure:product}.

We now give a way to compute the base of a factorized heap. If $H$ is a
heap, let the \emph{neighbourhood} of $H$, denoted $v(H)$, be the set of
positions concurrent to at least one piece of $H$.

\begin{lemma}\label{lemma:base}
Let $(H=H_1\cdot H_2)$ be a factorized heap. The base of the heap $H$ is given
by:
\[b(H)=b(H_1)\cup\bigl(b(H_2)\setminus v(H_1)\bigr).\]
\end{lemma}

\begin{proof}
As the heap $H$ is built by dropping all pieces of $H_1$, then all pieces of
$H_2$, no piece of $H_1$ can be above a piece of $H_2$. Therefore, all minimal
pieces of $H_1$ are also minimal pieces in $H$.

A minimal piece of $H_2$ is minimal in $H$ if and only if it is not above a
piece of $H_1$. This happens if and only if its position is not in the
neighbourhood of $H_1$, hence the given formula.
\end{proof}

Given a heap $H$ with base $S$, we may factorize $H$ as $H=S\cdot H_2$.
Lemma~\ref{lemma:base} asserts that the base of $H_2$ is included in $v(S)$.
This yields:
\begin{equation}\label{eq:neighbours}
\cH_S(t)=t^{\abs S}\cH_{[v(S)]}(t).
\end{equation}

\subsection{Heaps marked with a set of pieces}\label{subsection:marked}

A number of our problems in animal enumeration can be seen as enumeration of
heaps marked with a set of pieces; for example, computing the average number
of adjacent sites in animals is linked to enumerating animals with two
adjacent sites marked, which is in turn linked to enumerating heaps with some
pieces marked. Here, we give a means to link such marked heaps to factorized
heaps, which prove to be more manageable.

\begin{definition}\label{definition:marked}
A \emph{marked heap} $(H,X)$ is a heap $H$, marked with a set $X$ of
\emph{pairwise independent} pieces.
\end{definition}

\begin{definition}
Let $(H,X)$ be a marked heap. Let $\mathcal F_\downarrow(H,X)$ be the
factorized heap $(H=H_1\cdot H_2)$ where $H_1$ consists of all pieces below at
least a piece of $X$. Let $\mathcal F_\uparrow(H,X)$ be the factorized heap
$(H=H_1\cdot H_2)$ where $H_2$ consists of all pieces above at least a piece
of $X$.
\end{definition}

\begin{definition}\label{definition:almost-marked}
An \emph{almost strict marked heap} is a marked heap $(H,X)$ such that no
piece of $H$ sits on an unmarked piece at the same position.
\end{definition}

We call \emph{marked stack} a stack of an almost strict marked heap containing
a marked piece; such a stack may have one or two pieces, and the marked piece
is always the lowest piece of the stack.

\begin{definition}
Let $(H,X)$ be an almost strict marked heap. Let $X^+$ be the set consisting
of the highest piece of each marked stack. Define the following factorized
heaps:
\begin{align*}
F_\downarrow(H,X)&=\mathcal F_\downarrow(H,X)\text;\\
F_\uparrow(H,X)&=\mathcal F_\uparrow(H,X^+)\text.
\end{align*}
\end{definition}

\begin{figure}[ht]
\begin{center}
\begin{tikzpicture}[dominoes]
\begin{scope}[yshift=14cm]
\node at (0,0) [domino] {};
\node at (-1,1) [domino] {};
\node at (-1,2) [marked domino] {};
\node at (1,1) [domino] {};
\node at (0,3) [domino] {};
\node at (0,4) [domino] {};
\node at (4,0) [domino] {};
\node at (3,1) [domino] {};
\node at (2,2) [domino] {};
\node at (5,1) [domino] {};
\node at (5,2) [marked domino] {};
\node at (5,3) [domino] {};
\draw[arrow] (8,2) -- +(3,0);
\end{scope}
\begin{scope}[yshift=13cm,xshift=30cm]
\node at (0,0) [domino] {};
\node at (-1,1) [domino] {};
\node at (-1,4) [marked domino] {};
\node at (1,1) [domino] {};
\node at (0,5) [domino] {};
\node at (0,6) [domino] {};
\node at (4,0) [domino] {};
\node at (3,1) [domino] {};
\node at (2,2) [domino] {};
\node at (5,1) [domino] {};
\node at (5,4) [marked domino] {};
\node at (5,5) [domino] {};
\end{scope}
\begin{scope}[yshift=13cm,xshift=15cm]
\node at (0,0) [domino] {};
\node at (-1,1) [domino] {};
\node at (-1,2) [marked domino] {};
\node at (1,4) [domino] {};
\node at (0,5) [domino] {};
\node at (0,6) [domino] {};
\node at (4,0) [domino] {};
\node at (3,4) [domino] {};
\node at (2,5) [domino] {};
\node at (5,1) [domino] {};
\node at (5,2) [marked domino] {};
\node at (5,4) [domino] {};
\end{scope}
\begin{scope}[yshift=1cm]
%\node at (0,0) [domino] {};
%\node at (2,0) [domino] {};
\node at (-1,1) [domino] {};
\node at (1,1) [marked domino] {};
\node at (3,1) [domino] {};
\node at (1,2) [domino] {};
\node at (4,2) [domino] {};
\node at (0,3) [domino] {};
\node at (3,3) [marked domino] {};
\node at (5,3) [domino] {};
%\node at (-1,4) [domino] {};
%\node at (2,4) [domino] {};
\node at (4,4) [domino] {};
\node at (3,5) [domino] {};
\draw[arrow] (8,3) -- +(3,0);
\end{scope}
\begin{scope}[xshift=30cm]
%\node at (0,0) [domino] {};
%\node at (2,0) [domino] {};
\node at (-1,1) [domino] {};
\node at (1,1) [domino] {};
\node at (3,1) [domino] {};
\node at (1,5) [marked domino] {};
\node at (4,2) [domino] {};
\node at (0,6) [domino] {};
\node at (3,5) [marked domino] {};
\node at (5,3) [domino] {};
%\node at (-1,7) [domino] {};
%\node at (2,6) [domino] {};
\node at (4,6) [domino] {};
\node at (3,7) [domino] {};
\end{scope}
\begin{scope}[xshift=15cm]
%\node at (0,0) [domino] {};
%\node at (2,0) [domino] {};
\node at (-1,5) [domino] {};
\node at (1,1) [marked domino] {};
\node at (3,1) [domino] {};
\node at (1,5) [domino] {};
\node at (4,2) [domino] {};
\node at (0,6) [domino] {};
\node at (3,3) [marked domino] {};
\node at (5,5) [domino] {};
%\node at (-1,7) [domino] {};
%\node at (2,6) [domino] {};
\node at (4,6) [domino] {};
\node at (3,7) [domino] {};
\end{scope}
\end{tikzpicture}
\end{center}
\caption{A marked heap and an almost strict marked heap (left), their image by
$\mathcal F_\downarrow$ and $F_\downarrow$ (middle) and by $\mathcal
F_\uparrow$ and $F_\uparrow$ (right).}\label{figure:factor}
\end{figure}
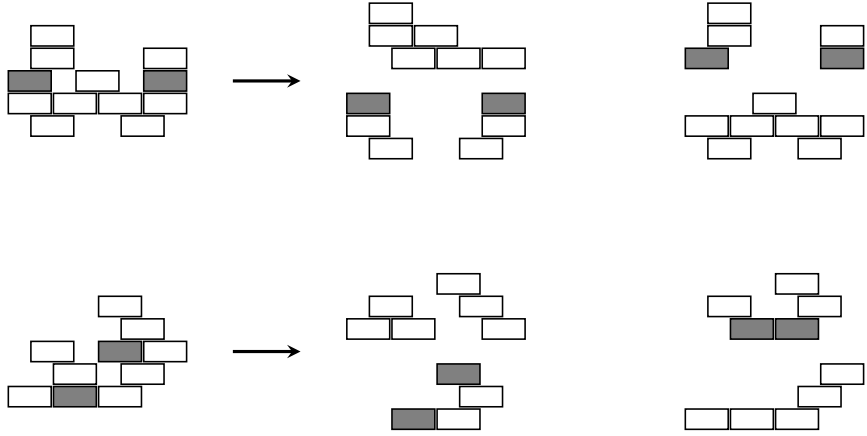

\begin{lemma}
The mappings $\mathcal F_\downarrow$ and $\mathcal F_\uparrow$ (resp.\
$F_\downarrow$ and $F_\uparrow$) are bijections from the set of marked heaps
to the set of factorized heaps (resp.\ almost strict marked heaps to almost
strict factorized heaps).

Their inverse bijections are as follows. Let $(H=H_1\cdot H_2)$ be a
factorized heap, let $X$ be the set of maximal pieces of $H_1$ and $Y$ the set
of minimal pieces of $H_2$. We have:
\[(H_1\cdot H_2)=\mathcal F_\downarrow(H,X)=\mathcal F_\uparrow(H,Y).\]
If $(H=H_1\cdot H_2)$ is an almost strict factorized heap, let $Y^-$ be the
set consisting of the lowest piece of each stack of $H$ containing a piece of
$Y$. We have:
\[(H_1\cdot H_2)=F_\downarrow(H,X)=F_\uparrow(H,Y^-).\]
\end{lemma}

\begin{proof}
Let us first do the case of marked heaps.
This fact easily stems from the poset structure on the pieces of a heap: the
set $X$ of maximal pieces of $H_1$ is the only set of pairwise independent
pieces such that every piece of $H_1$ is below a piece of $X$. The case of
$\mathcal F_\uparrow$ is identical.

We handle almost strict marked heaps similarly. Pulling downwards the lowest
piece of each marked stack, or pushing upwards the highest piece, ensures that
$H_1$ and $H_2$ are strict in the resulting factorized heap; conversely, as
the pieces of $X$ and $Y^-$ are the lowest of their stacks, $(H,X)$ and
$(H,Y^-)$ are almost strict marked heaps.
\end{proof}

As a first application, we give a means to compute the generating functions of
heaps marked with one piece at a fixed position.

\begin{definition}
Let $q$ be a position and $S$ a set of positions. We denote by
$\cH^{(q)}_{[S]}(t)$ the generating function of heaps with base included in
$S$, marked with a piece at position $q$. We denote by $\cV^q_{[S]}(t)$ the
set of heaps with base included in $S$ \emph{avoiding} $q$, i.e.\ such that no
piece is concurrent to $q$. As usual, we use analogous notations for strict
heaps.
\end{definition}

\begin{lemma}\label{lemma:marked}
The generating functions counting heaps of base included in $S$ marked with a
piece at position $q$ is given by:
\begin{align*}
\mathcal H^{(q)}_{[S]}(t)&=\begin{cases}
\mathcal H_{[S]}(t)\mathcal H_{\{q\}}(t)&\text{if }q\in S\text,\\
\bigl(\mathcal H_{[S]}(t)-\mathcal V^q_{[S]}(t)\bigr)
\mathcal H_{\{q\}}(t)
&\text{otherwise;}
\end{cases}\\
H^{(q)}_{[S]}(t)&=\frac1{1+t}\begin{cases}
H_{[S]}(t)H_{\{q\}}(t)&\text{if }q\in S\text,\\
\bigl(H_{[S]}(t)-V^q_{[S]}(t)\bigr)H_{\{q\}}(t)
&\text{otherwise.}
\end{cases}
\end{align*}
\end{lemma}

\begin{proof}
Let $(H,x)$ be a heap with a marked piece at position $q$. We use the
bijection $\mathcal F_\uparrow$ to turn it into a factorized heap $(H_1\cdot
H_2)$. We know that $H_2$ has base $\{q\}$, and that $H_1\cdot H_2$ has base
included in $S$.  According to Lemma~\ref{lemma:base}, we have: \[b(H_1\cdot
H_2)=b(H_1)\cup\bigl(\{q\}\setminus v(H_1)\bigr).\] If $q$ is in $S$, this
simply means that $b(H_1)$ is included in $S$; if not, it means that $q$ must
be in $v(H_1)$ as well, so that $H_1$ must not avoid $q$. We thus get the
result for general heaps.

Let $H^{(q)*}_{[S]}(t)$ be the generating function of almost strict marked
heaps, with exactly one marked piece at position $q$. The bijection $\mathcal
F_\uparrow$ turns these heaps into almost strict factorized heaps, on which we
apply the same reasoning. Moreover, as only one piece may be duplicated, we
have the identity:
\[H^{(q)*}_{[S]}(t)=(1+t)H^{(q)}_{[S]}(t).\]
This yields the second formula.
\end{proof}

\section{Directed animals and heaps of dominoes}\label{section:dominoes}

\subsection{Definitions}

\begin{definition}
Let $Q$ be either a subset of $\bZ$ or of the form $\bZ/m\bZ$ with $m$ an even
number. The \emph{square lattice} over $Q$, denoted by $\Gamma_Q$, is the
oriented graph with vertices $(q,i)\in Q\times\bN$ such that $q+i$ is even,
and edges $(q,i)\to(q+1,i+1)$ and $(q,i)\to(q-1,i+1)$. The \emph{triangular
lattice} over $Q$, denoted by $\Delta_Q$, is $\Gamma_Q$ with additionnal edges
$(q,i)\to(q,i+2)$.
\end{definition}

\begin{definition}
Let $Q$ be a subset or quotient of $\bZ$ in the same conditions as above. Let
$\cC$ be the relation defined by $(q,q')\in\cC$ if and only if
$\abs{q-q'}\leq1$. The model of heaps $(Q,\cC)$ is called the model of
\emph{heaps of dominoes} with set of positions $Q$.
\end{definition}

\begin{figure}[ht]
\begin{center}
\begin{tikzpicture}
\begin{scope}[scale=.3]
\draw[semithick,dashed] (-.5,0) -- ++(0,6.5) (5.5,0) -- ++(0,6.5);
\draw[arrow] (8,3) -- ++(3,0);
\begin{scope}
\clip (-.5,0) rectangle ++(6,6.5);
\draw[help lines,animal] (-3,-3) grid (7,7);
\end{scope}
\begin{scope}[animal,every node/.style={site}]
\node at (1,-1) {};
\node at (2,-1) {};
\node at (3,-1) {};
\node at (4,-1) {};
\node at (0,0) {};
\node at (1,0) {};
\node at (3,0) {};
\node at (2,2) {};
\node at (3,2) {};
\node at (3,3) {};
\end{scope}
\end{scope}
\begin{scope}[dominoes,xshift=14cm,every node/.style={domino}]
\draw[semithick,dashed] (-1.1,-.5) -- ++(0,7) (5.1,-.5) -- ++(0,7);
\begin{scope}
\clip (-1,-.5) rectangle ++(6,7);
\node at (0,0) {};
\node at (2,0) {};
\node at (1,1) {};
\node at (3,1) {};
\node at (4,2) {};
\node at (-1,3) {};
\node at (3,3) {};
\node at (5,3) {};
\node at (0,4) {};
\node at (1,5) {};
\node at (0,6) {};
\end{scope}
\end{scope}
\begin{scope}[scale=.3,yshift=-10cm]
\draw[arrow] (8,3) -- ++(3,0);
\foreach \x in {0,2,4}
{\draw[help lines] (\x,0) -- ++(0,6.5);}
\foreach \x in {1,3}
{\draw[help lines] (\x,1) -- ++(0,5.5);}
\begin{scope}
\clip (0,0) rectangle ++(4,6.5);
\draw[help lines,animal] (-3,-3) grid (7,7);
\end{scope}
\begin{scope}[animal,every node/.style={site}]
\node at (0,0) {};
\node at (1,0) {};
\node at (2,0) {};
\node at (3,0) {};
\node at (4,0) {};
\node at (2,1) {};
\node at (5,1) {};
\node at (3,2) {};
\node at (3,3) {};
\end{scope}
\end{scope}
\begin{scope}[dominoes,xshift=14cm,yshift=-10cm,every node/.style={domino}]
\clip (-1,-.5) rectangle ++(6,7);
\node at (0,0) {};
\node at (1,1) {};
\node at (2,2) {};
\node at (1,3) {};
\node at (3,3) {};
\node at (4,4) {};
\node at (1,5) {};
\node at (0,6) {};
\node at (4,6) {};
\end{scope}
\end{tikzpicture}
\end{center}
\caption{An animal on the square cylindrical lattice of width~$6$ and the
triangular rectangular lattice of width~$5$, and their corresponding heaps of
dominoes.}\label{figure:models}
\end{figure}
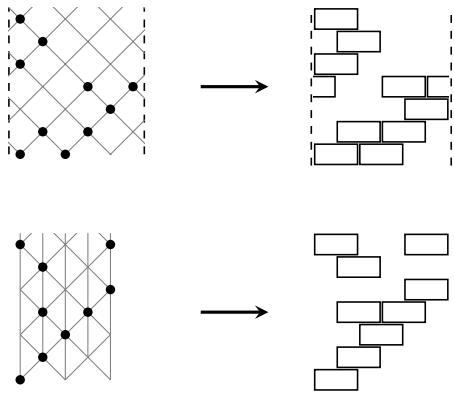

Up to a translation, there are four kinds of models and associated lattices:
\begin{itemize}
\item the \emph{full model} is the model $Q=\bZ$;
\item the \emph{half model} is the model $Q=\bN$;
\item the \emph{cylindrical model} of width $m$ is the model $Q=\bZ/m\bZ$,
with $m$ even;
\item the \emph{rectangular model} of width $m$ is the model
$Q=\{0,\dotsc,m-1\}$.
\end{itemize}
As seen in Figure~\ref{figure:models}, directed animals of source $S$ on the
square lattice $\Gamma_Q$ are identical to strict heaps of dominoes of base
$S$ in the model $(Q,\cC)$, while directed animals on the triangular lattice
$\Delta_Q$ are inflated heaps.

We say that a heap $H$ of dominoes is \emph{aligned} if all its pieces $(q,i)$
are such that $q+i$ is even. In particular, all heaps and inflated heaps
corresponding to directed animals are aligned.

Let $\Gamma_Q$ be a square lattice and $\Delta_Q$ its associated triangular
lattice; let $S$ be a one-line source, i.e., a set of vertices of the form
$(q,0)$. We denote by $A_S(t)$ and $\cA_S(t)$ the generating functions of
animals of source $S$ in the lattices $\Gamma_Q$ and $\Delta_Q$, respectively.
These two generating functions also count strict and inflated heaps of base
$S$ in the model of heaps of dominoes $Q$.

In this section, we give means to compute $\cA_S(t)$; the generating function
$A_S(t)$ is then given by performing the substitution
\smash{$t\mapsto\frac{t}{1+t}$}. Of course, we compute in the same way the
generating functions $\cA_{[S]}(t)$ and $A_{[S]}(t)$ of animals with a source
included in $S$.

\subsection{Bounded lattices}

The bounded lattices correspond to finite models of heaps of dominoes, that
is, the cylindrical models and the rectangular models. Let $Q$ be a finite
model and $S$ be a set of positions. The identity \eqref{eq:neighbours} and
Lemma~\ref{lemma:inversion} give the value of the generating function of heaps
of base $S$:
\[\cA_S(t)=t^{\abs S}\frac{\cT_{[Q\smallsetminus v(S)]}(t)}{\cT_{[Q]}(t)}.\]
All we need is therefore to compute the generating functions of trivial heaps.
Do do this, we define two sequences of polynomials.

\begin{definition}\label{definition:fm}
Define the sequences $(F_m(t))_{m\geq0}$ and $(\widehat F_m(t))_{m\geq2}$ of
polynomials by $F_0(t)=1$, $F_1(t)=1-t$, and for all $m\geq2$:
\begin{align*}
F_m(t)&=F_{m-1}(t)-tF_{m-2}(t)\text;\\
\widehat F_m(t)&=F_{m-1}(t)-tF_{m-3}(t)\text.
\end{align*}
\end{definition}

The polynomials $F_m(t)$ are often called the \emph{Fibonacci polynomials}.
With these polynomials, we can compute the generating function $\cT_{[S]}(t)$,
in any finite model $Q$ and for any set $S$, using the two lemmas below. We
state them without proof, and refer to \cite{viennot} for more detail.
Examples are also given in Section~\ref{section:examples}.

\begin{lemma}\label{lemma:trivial1}
Let $m\geq0$. The generating function $\cT_{[Q]}(t)$ of trivial heaps in the
rectangular model of width~$m$ is $F_m(t)$. If $m$ is even, the generating
function of trivial heaps in the cylindrical model of width~$m$ is $\widehat
F_m(t)$.
\end{lemma}

\begin{lemma}\label{lemma:trivial2}
Let $S$ be a finite set of positions; write $S=S_1\cup\dotsb\cup S_k$, where
the $S_i$ are intervals of $\mathbb Z$, with $k$ minimal. The generating
function of trivial heaps included in $S$ is:
\[\cT_{[S]}(t)=F_{\abs{S_1}}(t)\dotsm F_{\abs{S_k}}(t).\]
\end{lemma}

We now give explicit formul\ae{} for the special case $S=\{0\}$. In the
following, we call \emph{zero-source animals} the animals of source $\{0,0\}$
in any model. Let $\mathcal A_m(t)$ and $\mathcal D_m(t)$ be the generating
functions of zero-source animals in the cylindrical and rectangular triangular
lattices of width $m$, respectively. We have:
\begin{align}\label{am}
\mathcal A_m(t)&=\frac{F_{m-1}(t)}{\widehat F_m(t)}-1\text;\\\label{dm}
\mathcal D_m(t)&=\frac{F_{m-1}(t)}{F_m(t)}-1\text.
\end{align}
The generating functions $A_m(t)$ and $D_m(t)$ counting zero-source animals on
the square lattices are derived by performing the substitution
\smash{$t\mapsto\frac t{1+t}$}.

\subsection{Unbounded lattices}

We now address the unbounded lattices, i.e.\ the full and half lattices.
We start with zero-source animals, defined above, which are simplest. In the
rectangular and half models, such animals are also called \emph{half-animals}.

\begin{definition}
Let $A(t)$ and $D(t)$ be the generating functions of zero-source animals
on the full square lattice and the half square lattice, respectively. Let
$\cA(t)$ and $\cD(t)$ be their counterparts on the triangular lattices.
\end{definition}

These four generating functions are given below (see
\cite{dhar,gouyou,betrema}):

\begin{proposition}
The generating functions of zero-source animals on the infinite models are:
\begin{align}
\cA(t)&=\fulltriang\text;\\
\cD(t)&=\halftriang\text;\\
A(t)&=\fullsquare\text;\\
D(t)&=\halfsquare\text.
\end{align}
\end{proposition}
Of course, as usual, the generating functions counting animals on the square
lattices can be obtained by performing the substitution
$t\mapsto\frac{t}{1+t}$ in the ones counting animals on the triangular
lattices.

\begin{definition}\label{definition:compact}
In the full model, a \emph{compact source} is a finite set of consecutive even
positions. In the half model, a compact source is a finite set of consecutive
even positions, starting at 0.
\end{definition}

The next result gives the generating function of animals with a given compact
source. The proof may be found in \cite{betrema}.

\begin{proposition}\label{proposition:compact}
Let $C$ be a compact source with $k$ sites. The generating function of animals
of source $C$ on the full triangular lattice is:
\begin{align*}
\cA_C(t)&=\cD(t)^{k-1}\cA(t)\text.\\
\intertext{On the half lattice, this generating function is:}
\cA_C(t)&=\cD(t)^k\text.
\end{align*}
\end{proposition}

\begin{remark}
From this, it can be proved that the number of animals of area~$n$ with
\emph{any} compact source on the triangular lattice is $4^{n-1}$, and
$3^{n-1}$ on the square lattice. See \cite{gouyou,betrema}.
\end{remark}

Finally, we are able to compute the generating function of animals with an
arbitrary source $S$.

\begin{proposition}
Let $Q=\bZ$ or $\bN$. Let $S\subseteq Q$ be a one-line source and $C$ be the
smallest compact source containing $S$. The generating function of animals of
source $S$ in $\Delta_Q$ is:
\[\cA_S(t)=t^{\abs S-\abs C}\cT_{[v(C)\smallsetminus v(S)]}(t)\cA_C(t).\]
\end{proposition}

\begin{proof}
Let us address the full lattice. Let $Q_m$ be a cylindrical model large enough
so that $C\subseteq Q_m$. The generating function of heaps of base $S$ in the
model $Q_m$ is:
\[\cA_{S,m}(t)=
t^{\abs S}\frac{\cT_{[Q_m\smallsetminus v(S)]}(t)}{\cT_{[Q_m]}(t)}.\]
As $C$ is the smallest compact source containing $S$, $v(C)$ is also the
smallest interval containing $v(S)$, which ensures that no position of
$Q_m\setminus v(C)$ can be concurrent to $v(C)\setminus v(S)$. Therefore,
Lemma~\ref{lemma:trivial2} entails that:
\[\cT_{[Q_m\smallsetminus v(S)]}(t)=
\cT_{[v(C)\smallsetminus v(S)]}(t)\cT_{[Q_m\smallsetminus v(C)]}(t),\]
and thus:
\[\cA_{S,m}(t)=
t^{\abs S-\abs C}\cT_{[v(C)\smallsetminus v(S)]}(t)\cA_{C,m}(t).\]
We conclude by letting $m$ tend to infinity.

In the case of the half-lattice, we repeat the same reasoning, this time
taking for $Q_m$ a rectangular model large enough so that $C\subseteq Q_m$.
\end{proof}

\section{Average number of adjacent sites, loops and neighbours of directed
animals}\label{section:results}

\subsection{Notations and results}

In this section, $\Gamma$ is a square lattice (full, half, cylindrical or
rectangular), and $S$ is a one-line source of $\Gamma$. In
Section~\ref{section:dominoes}, we have computed the generating functions
$A_S(t)$ and $A_{[S]}(t)$ counting directed animals on $\Gamma$ of source $S$
and with a source included in $S$.

We are now interested in three parameters of the directed animals: number of
adjacent sites, number of loops, and site perimeter. The first two are
defined in Section~1; we denote by $j(A)$ and $\ell(A)$ the number of pairs of
adjacent sites and number of loops of the animal $A$, respectively.

\renewcommand\S{$S$\nobreakdash}

\begin{definition}\label{perim}
Let $A$ be an animal on $\Gamma$ with a source included in $S$. An
\emph{\S-neighbour} of $A$ is a vertex $v$ of $\Gamma$ not in $A$
such that $A\cup\{v\}$ is still a directed animal with a source included in
$S$.

Assume now that the graph $\Gamma$ is embedded in a larger graph $\Gamma'$. An
\emph{internal \S-neighbour} of $A$ is an \S-neighbour of $A$ seen as an
animal on $\Gamma$. An \emph{external \S-neighbour} of $A$ is an
\S-neighbour of $A$ seen as an animal on $\Gamma'$.
% When there is no ambiguity on the set $S$, we abusively write \emph{neighbour} istead of \S-neighbour.

Finally, a vertex of $\Gamma$ is at the \emph{edge} of the lattice if it has
outdegree~1.
\end{definition}

For the purpose of this definition, we regard the half-lattice and the
rectangular lattices as embedded in the full lattice. The full and cylindrical
lattices are not naturally embedded in any larger graph, so in these lattices
internal and external neighbours are identical. Moreover, if $S$ is the source
of $A$, the \S-neighbours of $A$ coincide with the usual neighbours of $A$, as
defined in Section~1.

Assuming no ambiguity on the set $S$, we denote by $p_i(A)$ the number of
internal \S-neighbours of $A$ (or internal site perimeter) and by $p_e(A)$ its
number of external \S-neighbours (or external site perimeter). We also denote
by $e(A)$ the number of sites of $A$ at the edge of the lattice $\Gamma$.

The generating functions defined below are linked to the average value of each
parameter in animals of a given area.

\begin{definition}\label{definition:square1}
Define the following generating functions, counting animals with a source
included in $S$:
\begin{itemize}
\item Let $J_{[S]}(t)$ be the generating function of animals marked with two
adjacent sites.
\item Let $L_{[S]}(t)$ be the generating function of animals marked with a
loop.
\item Let $P^e_{[S]}(t)$ be the generating function of animals marked with an
external \S-neighbour.
\item Let $P^i_{[S]}(t)$ be the generating function of animals marked with an
internal \S-neighbour.
\end{itemize}
Also define $J_S(t)$, $L_S(t)$, $P^e_S(t)$ and $P^i_S(t)$ the analogous
generating functions of animals of source $S$.
\end{definition}

To compute these generating functions, we again use the correspondence between
animals and heaps of dominoes. We denote by $(Q,\mathcal C)$ the model of
heaps of dominoes associated with the lattice $\Gamma$; we also denote by $S$
the aligned set of positions associated to the source $S$.

We now define some generating functions counting heaps. As usual, a generating
function with a subscript $[S]$ counts heaps with a base included in $S$, and
one with a subscript $S$ counts heaps with base $S$, so this precision will
often be omitted from the definition.

In some cases, we consider heaps having a minimal piece outside $S$; we call
such a piece an \emph{illegal minimal piece}. Note that an illegal minimal
piece does not have to be aligned with $S$.

\begin{definition}\label{definition:square2}
Define the following generating functions:
\begin{itemize}
\item Let $M_{[S]}(t)$ be the generating function of strict heaps marked with
a piece at a position $q$, such that $q+2$ is in $Q$.
\item If $q$ is in $S$, let $W^q_{[S]}(t)$ be the generating function of
strict heaps with a minimal piece at position $q$ and an illegal minimal piece
at $q+2$. Let $W_{[S]}(t)$ be the sum of $W^q_{[S]}(t)$ over all $q\in S$ such
that $q+2\not\in S$.
\item Let $E_{[S]}(t)$ be the generating function of strict heaps marked with
a piece at a position $q$, such that either $q-1$ or $q+1$ is not in $Q$.
\end{itemize}
Also define the analogues $M_S(t)$, $W_S(t)$ and $E_S(t)$; thus, $W_S(t)$ is
the generating functions of strict heaps of base $S\cup\{q+2\}$, such that $q$
is in $S$ but $q+2$ is not.
\end{definition}

We now show that these generating functions can be computed using previous
results. First, we have:
\begin{equation}\label{eq:Wcalc}
W_{[S]}(t)=\sum_{\substack{q\in S\\q+2\not\in S}}
\sum_{\substack{T\subseteq S\\q\in T}}A_{T\cup\{q+2\}}(t),
\end{equation}
where $A_{T\cup\{q+2\}}(t)$ counts strict heaps with base $T\cup\{q+2\}$, and
is computed using the results of Section~3. Next, the generating function of
strict heaps with base included in $S$ marked with a single piece is
$tA_{[S]}'(t)$. Thus, we have:
\begin{align}\label{eq:Mcalc}
M_{[S]}(t)&=tA_{[S]}'(t)-\sum_{q+2\not\in Q}A_{[S]}^{(q)}(t)\text;\\
\label{eq:Ecalc}
E_{[S]}(t)&=\sum_{q-1\not\in Q\text{ or }q+1\not\in Q}A_{[S]}^{(q)}(t)\text,
\end{align}
where $A_{[S]}^{(q)}(t)$ counts strict heaps with base included in $S$ marked
with a piece at position $q$ and is computed using Lemma~\ref{lemma:marked}.
In all three equations, the sum goes over a finite number of positions $q$,
which ensures that all three generating functions can be computed. In this
regard, the full and cylindrical models are the simplest, as $M_{[S]}(t)$ is
equal to $tA_{[S]}'(t)$ and $E_{[S]}(t)$ is zero.

\begin{theorem}\label{theorem:square}
In square lattices, the generating functions counting the total number of
adjacent pieces, loops and site perimeters of the animals with source included
in $S$ are given by:
\begin{align}\label{Js}
J_{[S]}(t)&=\frac{tM_{[S]}(t)-W_{[S]}(t)}{1+t}\text;\\\label{Ls}
L_{[S]}(t)&=t(1+t)J_{[S]}(t)\text;\\\label{Pes}
P^e_{[S]}(t)&=\abs SA_{[S]}(t)+tA_{[S]}'(t)-J_{[S]}(t)\text;\\\label{Pis}
P^i_{[S]}(t)&=P^e_{[S]}(t)-E_{[S]}(t)\text.
\end{align}
Moreover, the corresponding generating functions for animals of source $S$ are:
\begin{align}\label{JS}
J_S(t)&=\frac{tM_S(t)+j(S)A_S(t)-W_S(t)}{1+t}\text;\\\label{LS}
L_S(t)&=t(1+t)J_S(t)\text;\\\label{PeS}
P^e_S(t)&=\abs SA_S(t)+tA_S'(t)-J_S(t)\text;\\\label{PiS}
P^i_S(t)&=P^e_S(t)-E_S(t)\text,
\end{align}
where $j(S)$ denotes the number of pairs of adjacent sites in the source $S$.
\end{theorem}

Applications of this theorem are given in Section~\ref{section:examples}.

\subsection{Site perimeters}\label{section:perimeter}

We first prove the four identities \eqref{Pes}, \eqref{Pis}, \eqref{PeS} and
\eqref{PiS} dealing with the internal and external site perimeter.

First, we remark that a vertex $(q,i)$ of $\Gamma$ has outdegree~1 if
and only if either $q-1$ or $q+1$ is not in $Q$. Thus, the generating function
$E_[S](t)$ satisfies:
\begin{equation}\label{Es}
E_{[S]}(t)=\sum_Ae(A)t^{\abs A}\text,
\end{equation}
where the sum goes over all animals of source included in $S$. The same goes
for the generating function $E_S(t)$.

\begin{lemma}\label{lemma:perimeter}
The number of external and internal \S-neighbours of every directed animal $A$
with source included in $S$ satisfy:
\begin{align*}
p_e(A)&=\abs S+\abs A-j(A)\text;\\
p_i(A)&=p_e(A)-e(A)\text.
\end{align*}
\end{lemma}

By summing the identities of this lemma over all animals of source included in
$S$, and using \eqref{Es}, we prove the identities \eqref{Pes} and
\eqref{Pis}. By summing them over all animals of source $S$, we get
\eqref{PeS} and \eqref{PiS}.

\begin{proof}
When dealing with the external site perimeter, the lattice $\Gamma$ is
embedded in a lattice $\Gamma'$ which is either the full lattice or a
cylindrical lattice. Let $Z$ be the number of pairs of vertices $(v,w)$ such
that $v$ is a site of $A$ and $w$ is a \emph{child} of $v$ (i.e., $v\to w$ is
an edge of $\Gamma'$), whether in $A$ or not. As every vertex has outdegree~2,
we have
\[Z=2\abs A.\]
Now, as $A$ is a directed animal, a child of a site of $A$ is either a site of
$A$ or an external \S-neighbour of $A$. The only sites and \S-neighbours of
$A$ not counted are the ones in $S$; moreover, two sites have a child in
common if and only if they are adjacent. Hence:
\[Z=\abs A+p_e(A)+j(A)-\abs S\text,\]
which yields the announced formula for $p_e(A)$.

If $\Gamma$ is either the half-lattice or a rectangular lattice, then each
site on the edge of the lattice has one external neighbour not in $\Gamma$.
Thus, we have:
\[p_e(A)=p_i(A)+e(A).\qedhere\]
\end{proof}

\subsection{Average number of adjacent sites and loops}\label{subsection:proof}

To prove the remaining identities of Theorem~\ref{theorem:square}, dealing
with the number of adjacent sites and loops, we use several bijections between
various sets of heaps marked with certain pieces.  Rather than strict marked
heaps, it is convenient here to use almost strict marked heaps (see
Definition~\ref{definition:almost-marked}).

Although the proof seems complicated due to the high number of sets and
associated generating functions that we must consider, each bijection is
actually very simple, consisting in adding/removing a single piece.

\begin{figure}[ht]
\begin{center}
\begin{tikzpicture}[dominoes]
\begin{scope}
\node at (-1,0) [domino] {};
\node at (1,0) [domino] {};
\node at (0,1) [domino] {$x$};
\draw [help lines] (-2,-.5) -- ++(-2,-2);
\draw [help lines] (2,-.5) -- ++(2,-2);
\node at (0,-2) {$H_1$};
\draw [arrow] (0,-4) -- ++(0,-3);
\end{scope}
\begin{scope}[yshift=-11cm]
\node at (-1,0) [domino] {$y$};
\node at (1,0) [domino] {$z$};
\node at (0,1) [domino,help lines] {};
\draw [help lines] (-2,-.5) -- ++(-2,-2);
\draw [help lines] (2,-.5) -- ++(2,-2);
\node at (0,-2) {$H_1'$};
\end{scope}
\begin{scope}[xshift=11cm]
\node at (-1,2) [domino] {$x$};
\node at (1,0) [domino] {$y$};
\draw [help lines] (-2,1.5) -- ++(-2,-2);
\draw [help lines] (2,-.5) -- ++(2,-2);
\draw [help lines] (0,1.5) -- ++(0,-2);
\node at (0,-2) {$H_1$};
\draw [arrow] (0,-4) -- ++(0,-3);
\end{scope}
\begin{scope}[xshift=11cm,yshift=-11cm]
\node at (-2,1) [domino] {$z$};
\node at (-1,2) [domino, help lines] {};
\node at (1,0) [domino] {$y$};
\draw [help lines] (-3,.5) -- ++(-2,-2);
\draw [help lines] (2,-.5) -- ++(2,-2);
\draw [help lines] (-1,.5) -- ++(1,-1);
\node at (-.5,-2) {$H_1'$};
\end{scope}
\begin{scope}[xshift=21cm]
\node at (-1,0) [domino] {$y$};
\node at (1,2) [domino] {$x$};
\draw [help lines] (-2,-.5) -- ++(-2,-2);
\draw [help lines] (2,1.5) -- ++(2,-2);
\draw [help lines] (0,1.5) -- ++(0,-2);
\node at (0,-2) {$H_1$};
\draw [arrow] (0,-4) -- ++(0,-3);
\end{scope}
\begin{scope}[xshift=21cm,yshift=-11cm]
\node at (-1,0) [domino] {$y$};
\node at (1,2) [domino, help lines] {};
\node at (2,1) [domino] {$z$};
\draw [help lines] (-2,-.5) -- ++(-2,-2);
\draw [help lines] (3,.5) -- ++(2,-2);
\draw [help lines] (1,.5) -- ++(-1,-1);
\node at (.5,-2) {$H_1'$};
\end{scope}
\begin{scope}[xshift=32cm,yshift=-2cm]
\node at (0,4) [domino] {$x$};
\node at (2,0) [marked domino] {$y$};
\draw [help lines] (-1,3.5) -- ++(-2,-2);
\draw [help lines] (1,3.5) -- ++(0,-4);
\draw [help lines] (-3,-.5) -- ++(7,0);
\node at (-1,1) {$H_1$};
\draw [arrow] (0,-2) -- ++(0,-3);
\end{scope}
\begin{scope}[xshift=32cm,yshift=-13cm]
\node at (-1,3) [domino] {$z$};
\node at (0,4) [domino, help lines] {};
\node at (2,0) [marked domino] {$y$};
\draw [help lines] (-2,2.5) -- ++(-2,-2);
\draw [help lines] (0,2.5) -- ++(1,-1) -- ++(0,-2);
\draw [help lines] (-4,-.5) -- ++(8,0);
\node at (-1,1) {$H_1'$};
\end{scope}
\draw [brace] (7,3.5) -- ++(18,0);
\node at (0,4) {$\bL^*_{[S]}$};
\node at (16,5) {$\bI^{2*}_{[S]}\setminus\bJ^*_{[S]}$};
\node at (32,4) {$\bX^{2*}_{[S]}\setminus\bW^*_{[S]}$};
\draw [brace] (26,-14.5) -- ++(-20,0);
\node at (0,-15) {$\bJ^*_{[S]}$};
\node at (16,-16) {$\bI^{3*}_{[S]}$};
\node at (32,-15) {$\bX^{3*}_{[S]}$};
\end{tikzpicture}
\end{center}
\caption{On the left, the bijection $\Phi_0$: removing the piece $x$ uncovers
two adjacent pieces $y$ and $z$.
In the middle, the bijection $\Phi_1$: removing the piece $x$ uncovers a piece
$z$, with position at distance~3 from $y$ and higher than $y$.
On the right, the bijection $\Phi_2$, with the illegal minimal piece $y$
colored gray. Removing the piece $x$ uncovers a piece $z$, with position at
distance~3 from $y$.}\label{figure:bijection1}
\end{figure}
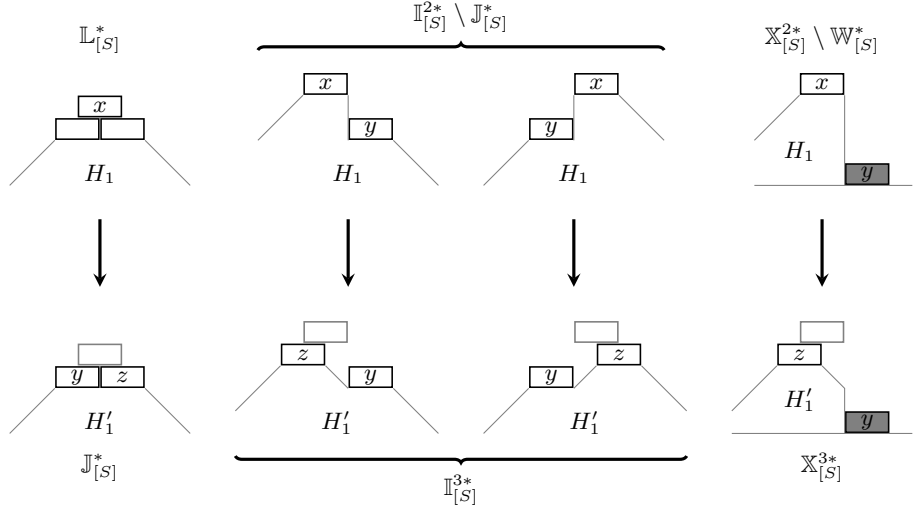

\begin{figure}[ht]
\newcommand\x{$x\scriptstyle+$}
\newcommand\y{$y\scriptstyle+$}
\begin{center}
\begin{tikzpicture}[dominoes]
\begin{scope}
\node at (-1,0) [domino] {\x};
\node at (1,0) [domino] {\y};
\draw [help lines] (-2,.5) -- ++(-2,2);
\draw [help lines] (2,.5) -- ++(2,2);
\node at (0,2) {$H_2$};
\node at (0,5) {$\bI^{2*}_{[S]}$};
\end{scope}
\begin{scope}[xshift=11cm]
\node at (-1,0) [domino] {\x};
\node at (1,0) [marked domino] {\y};
\draw [help lines] (-2,.5) -- ++(-2,2);
\draw [help lines] (2,.5) -- ++(2,2);
\node at (0,2) {$H_2$};
\node at (0,5) {$\bX^{2*}_{[S]}$};
\end{scope}
\begin{scope}[xshift=-5cm,yshift=-10cm]
\node at (0,0) [domino] {\x};
\node at (2,0) [domino, help lines] {};
\draw [help lines] (-1,.5) -- ++(-2,2);
\draw [help lines] (1,.5) -- ++(3,3);
\node at (.5,3) {$H_2'$};
\node at (.5,-2) {$\bM^*_{[S]}$};
\end{scope}
\begin{scope}[xshift=5cm,yshift=-10cm]
\node at (-1,0) [domino] {\x};
\node at (1,0) [domino, help lines] {};
\node at (2,1) [domino] {$z$};
\draw [help lines] (-2,.5) -- ++(-2,2);
\draw [help lines] (3,1.5) -- ++(2,2);
\node at (.5,3) {$H_2'$};
\node at (.5,-2) {$\bI^{3*}_{[S]}$};
\end{scope}
\begin{scope}[xshift=16cm,yshift=-10cm]
\node at (-1,0) [domino] {\x};
\node at (1,0) [domino, help lines] {};
\node at (2,1) [marked domino] {$z$};
\draw [help lines] (-2,.5) -- ++(-2,2);
\draw [help lines] (3,1.5) -- ++(2,2);
\node at (.5,3) {$H_2'$};
\node at (.5,-2) {$\bX^{3*}_{[S]}$};
\end{scope}
\draw [arrow] (5.5,-2) -- ++(0,-3);
\end{tikzpicture}
\end{center}
\caption{The bijection $\Phi_3$; illegal minimal pieces are colored gray.
After removing the piece $y^+$, either a new minimal piece $z$ of $H_2$ is
uncovered at position 3 or not. If it is, the heap $H'$ may be in
$\bI^{3*}_{[S]}$ or $\bX^{3*}_{[S]}$, depending on whether $z$ is an illegal
minimal piece. If not, the heap $H'$ is in $\bM^*_{[S]}$.}
\label{figure:bijection2}
\end{figure}
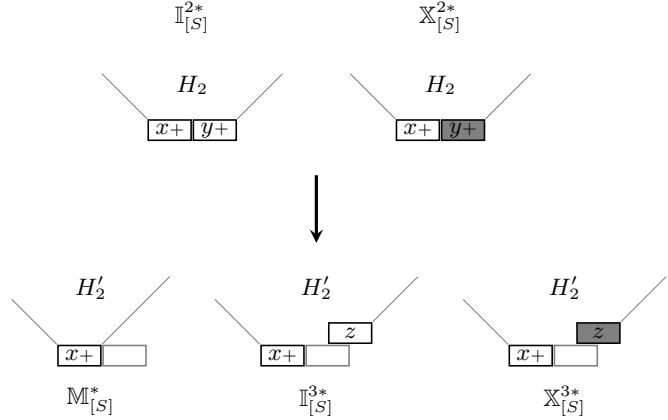

Each of these sets of heaps is illustrated in Figures \ref{figure:bijection1}
and \ref{figure:bijection2} with the relevant bijections.

\begin{definition}\label{definition:square3}
Define the following sets of almost strict marked heaps, assumed to have a
base included in $S$:
\begin{itemize}
\item Let $\bJ^*_{[S]}$ be the set of almost strict heaps marked with two
adjacent pieces.
\item Let $\bL^*_{[S]}$ be the set of almost strict heaps marked with the top
piece of a loop.
\item Let $\bM^*_{[S]}$ be the set of almost strict heaps marked with a piece
at a position $q$, such that $q+2$ is in $Q$.
\item Let $\bW^*_{[S]}$ be the set of almost strict heaps marked with a
minimal piece at a position $q$, and having an illegal minimal piece at
position $q+2$.
\item Let $\bI^{2*}_{[S]}$ (resp.\ $\bI^{3*}_{[S]}$) be the set of almost
strict heaps marked with two independent pieces at positions $q$ and $q+2$
(resp.\ $q+3$).
\item Let $\bX^{2*}_{[S]}$ (resp.\ $\bX^{3*}_{[S]}$) be the set of almost
strict heaps marked with a piece $x$ at a position $q$, and having an illegal
minimal piece at position $q+2$ (resp.\ $q+3$) independent from $x$.
\end{itemize}
Let $J^*_{[S]}(t)$, $L^*_{[S]}(t)$, $M^*_{[S]}(t)$, $W^*_{[S]}(t)$,
$I^{2*}_{[S]}(t)$, $I^{3*}_{[S]}$, $X^{2*}_{[S]}$, $X^{3*}_{[S]}(t)$ be the
generating functions of these sets.
\end{definition}

The asterisks $*$ are used above to denote sets of almost strict marked heaps;
we link these heaps to strict marked heaps later. Note that we have the
inclusions $\bJ^*_{[S]}\subseteq\bI^{2*}_{[S]}$ and
$\bW^*_{[S]}\subseteq\bX^{2*}_{[S]}$.
%Moreover, as the generating functions of
%almost strict marked heaps differ from the ones of strict marked heaps by a
%$1+t$ factor for each marked piece, we have:

\begin{lemma}\label{lemma:phi0}
The following identity holds:
\[L^*_{[S]}(t)=tJ^*_{[S]}(t).\]
\end{lemma}

\begin{proof}
To prove this result, we use a first bijection, which removes exactly one
piece:
\[\Phi_0\colon\bL^*_{[S]}\to\bJ^*_{[S]}.\]
Let $(H,\{x\})$ be a heap of $\bL^*_{[S]}$. Pull the marked piece $x$
downwards to form the factorized heap $(H_1\cdot H_2)=F_\downarrow(H,\{x\})$,
such that $x$ is the only maximal piece of $H_1$. Let $H_1=H_1'\cdot x$ and
$H'=H_1'\cdot H_2$. As $x$ is a loop, the heap $H_1'$ has two maximal pieces
$y$ and $z$, which are adjacent (see Figure~\ref{figure:bijection1}, left).
Moreover, as $H_1$ and $H_1'$ have same base and neighbourhood,
Lemma~\ref{lemma:base} guarantees that $H$ and $H'$ have the same base. We may
thus set:
\[\Phi_0\bigl(H,\{x\}\bigr)=\bigl(H',\{y,z\}\bigr).\]
This operation is easily reversible, by putting back the piece $x$.
\end{proof}

\begin{lemma}\label{lemma:phi1}
The following identity holds:
\[I^{2*}_{[S]}(t)-J^*_{[S]}(t)=tI^{3*}_{[S]}(t).\]
\end{lemma}

\begin{proof}
We again prove this result with a bijection removing one piece:
\[\Phi_1\colon\bI^{2*}_{[S]}\setminus\bJ^*_{[S]}\to\bI^{3*}_{[S]},\]
Let $(H,\{x,y\})$ be a heap of $\bI^{2*}_{[S]}\setminus\bJ^*_{[S]}$. We set
$(H_1\cdot H_2)=F_\downarrow(H,\{x,y\})$, pulling the pieces $x$ and $y$
downwards.
 
As the pieces $x$ and $y$ are not adjacent, one of them (say, $x$) is higher
than the other. Let $H_1=H_1'\cdot x$. As $H_1$ is strict, $H_1'$ must have a
second maximal piece $z$, such that the positions of $y$ and $z$ are at
distance~3 (Figure~\ref{figure:bijection1}, middle). Let $H'=H_1'\cdot H_2$.
As $H_1$ and $H_1'$ have same base and neighbourhood, we may set:
\[\Phi_1(H,\{x,y\})=(H',\{y,z\}).\]
This operation is reversible: let $(H',\{y,z\})$ be a heap of
$\bI^{3*}_{[S]}$, and let $(H_1'\cdot H_2)=F_\downarrow(H',\{y,z\})$. As the
set of positions $S$ is aligned, the heap $H_1'$ is also aligned. Therefore,
$y$ and $z$ cannot be at the same height; say, $z$ is higher. We set
$H_1=H_1'\cdot x$ so that $x$ sits on $z$ and the positions of $x$ and $y$ are
at distance~2, and $H=H_1\cdot H_2$; thus, we have
$\Phi_1(H,\{x,y\})=(H',\{y,z\})$.
\end{proof}

\begin{lemma}\label{lemma:phi2}
The following identity holds:
\[X^{2*}_{[S]}(t)-W^*_{[S]}(t)=tX^{3*}_{[S]}(t).\]
\end{lemma}

The generating function $W_T(t)$ may be found in
Definition~\ref{definition:square2}.

\begin{proof}
We use another bijection removing one piece:
\[\Phi_2\colon\bX^{2*}_{[S]}\setminus\bW^*_{[S]}\to\bX^{3*}_{[S]},\]
Let $(H,\{x\})$ be a heap of $\bX^{2*}_{[S]}\setminus\bW^*_{[S]}$, and let $y$
be the illegal minimal piece of $H$. We pull the pieces $x$ and $y$ downwards,
forming the factorized heap $(H_1\cdot H_2)$. As $x$ is not a minimal piece of
$H_1$, it must sit on another piece $z$, at position $q-1$
(Figure~\ref{figure:bijection1}, right). Let $H_1=H_1'\cdot x$ and
$H'=H_1'\cdot H_2$. Again, $H_1$ and $H_1'$ have the same base and
neighbourhood, so that $y$ is still the only illegal minimal piece of $H'$. We
set:
\[\Phi_2\bigl(H,\{x\}\bigr)=\bigl(H',\{z\}\bigr).\]
This operation is easily reversible by putting back the piece $x$.
\end{proof}

\begin{lemma}\label{lemma:phi3}
The following identity holds:
\[I^{2*}_{[S]}(t)+X^{2*}_{[S]}(t) =
t\Bigl(M^*_{[S]}(t)+I^{3*}_{[S]}(t)+X^{3*}_{[S]}(t)\Bigr).\]
\end{lemma}

\begin{proof}
We use a fourth and final bijection removing one piece:
\[\Phi_3\colon\bI^{2*}_{[S]}\cup\bX^{2*}_{[S]}
\to\bM^*_{[S]}\cup\bI^{3*}_{[S]}\cup\bX^{3*}_{[S]}.\]

In this proof, if $x$ is a piece, we write $x^+$ to denote the highest piece
of the stack of $x$, and $x^-$ to denote the lowest piece in the stack of $x$
(see Section~\ref{subsection:marked}).

Let $H$ be a heap of $\bI^{2*}_{[S]}$ or $\bX^{2*}_{[S]}$. In the first case,
let $x$ be the left-hand marked piece, and $y$ the right-hand one. In the
second, let $x$ be the marked piece and $y$ the illegal maximal piece. In
both cases, we set $(H_1,H_2)=F_\uparrow(H,\{x,y\})$, pushing upwards the
pieces $x^+$ and $y^+$.

Now, let $H_2=y^+\cdot H_2'$ and $H'=H_1\cdot H_2'$. We distinguish three
cases, illustrated in Figure~\ref{figure:bijection2}:
\begin{enumerate}[label=(\alph*)]
\item The piece $x^+$ is the only minimal piece of $H_2'$: the heap
$(H',\{x\})$ is in $\bM^*_{[S]}$.
\item The heap $H_2'$ has a minimal piece $z$ at position $q+3$, and $z$ is
not an illegal minimal piece of $H'$: the heap $(H',\{x,z^-\})$ is in
$\bI^{3*}_{[S]}$.
\item The heap $H_2'$ has a minimal piece $z$ at position $q+3$, and $z$ is an
illegal minimal piece of $H'$: the heap $(H',\{x\})$ is in $\bX^{3*}_{[S]}$.
\end{enumerate}
Once again, this operation is easily reversible by putting back the piece
$y^+$ and checking whether it is an illegal minimal piece.
\end{proof}

Finally, to prove the identity \eqref{JS}, we need a last lemma, given below.

\begin{lemma}\label{lemma:inclusion}
The following identity holds:
\[W_{[S]}(t)=\sum_{T\subseteq S}\bigl(W_T(t)-j(T)A_T(t)\bigr).\]
\end{lemma}

\begin{proof}
Consider the generating function:
\[\sum_{T\subseteq S}W_T(t).\]
We write $W_T(t)$ as the sum of all $W^q_T(t)$, for $q\in T$ and
$q+2\not\in T$. We then split the sum according to whether $q+2$ is in $S$ or
not:
\[\sum_{T\subseteq S}W_T(t)=\sum_{T\subseteq S}
\Biggl(\sum_{\substack{q\in T\\q+2\notin S}}W^q_T(t)\Biggr)+
\sum_{T\subseteq S}
\Biggl(\sum_{\substack{q\in T\\q+2\notin T\\q+2\in S}}W^q_T(t)\Biggr).\]
In the first term of the right-hand side of this equation, we recognize the
generating function $W_{[S]}(t)$. We rewrite the second term using the fact
that $W^q_T(t)=A_{T\cup\{q+2\}}(t)$ and by posing $T'=T\cup\{q+2\}$:
\begin{align*}
\sum_{T\subseteq S}W_T(t)&=W_{[S]}(t)+\sum_{T'\subseteq S}
\Biggl(\sum_{\substack{q\in T'\\q+2\in T'}}A_{T'}(t)\Biggr);\\
&=W_{[S]}(t)+\sum_{T'\subseteq S}j(T')A_{T'}(t).
\end{align*}
The lemma follows.
\end{proof}

With the above lemmas, we are now able to prove the theorem.

\begin{proof}[Proof of Theorem~\ref{theorem:square}]
The identities \eqref{Pes}, \eqref{Pis}, \eqref{PeS} and \eqref{PiS} dealing
with the site perimeters are proved in Section~\ref{section:perimeter}.

To prove the remaining identities, we first link the generating function
counting almost strict marked heaps with the ones counting strict marked
heaps. As each marked piece accounts for a $1+t$ factor, we have:
\begin{align*}
J^*_{[S]}(t)&=(1+t)^2J_{[S]}(t);\\
L^*_{[S]}(t)&=(1+t)L_{[S]}(t);\\
M^*_{[S]}(t)&=(1+t)M_{[S]}(t);\\
W^*_{[S]}(t)&=(1+t)W_{[S]}(t).
\end{align*}
Moreover, putting together the identities of Lemmas \ref{lemma:phi0},
\ref{lemma:phi1}, \ref{lemma:phi2} and \ref{lemma:phi3}, we find:
\begin{align*}
L^*_{[S]}(t)&=tJ^*_{[S]}(t);\\
J^*_{[S]}(t)&=tM^*_{[S]}(t)-W^*_{[S]}(t).
\end{align*}
Thus, we derive the first two identities of the theorem, \eqref{Js} and
\eqref{Ls}.

To prove the identities dealing with animals of source $S$, we remark that the
generating function $J_{[S]}(t)$ verifies:
\[J_{[S]}(t)=\sum_{T\subseteq S}J_T(t).\]
The generating functions $L_{[S]}(t)$ and $M_{[S]}(t)$ also behave in this
manner. Using the inclusion-exclusion principle, this means that the equation
\eqref{LS} giving $L_S(t)$ is a consequence of \eqref{Ls}.

To address the generating function $W_{[S]}(t)$, we use
Lemma~\ref{lemma:inclusion}, rewriting \eqref{Js} as:
\[\sum_{T\subseteq
S}\biggl(J_T(t)-\frac{tM_T(t)+j(T)A_T(t)-W_T(t)}{1+t}\biggr) = 0.\]
The identity \eqref{JS} is then derived using the inclusion-exclusion
principle.
\end{proof}

\subsection{Triangular lattices}

Let $\Delta$ be the triangular lattice corresponding to $\Gamma$. A number of
our results on the animals of $\Gamma$ have counterparts on the animals of
$\Delta$. The results and proofs are very similar, and we go into slightly
less detail.

Given an animal $A$ on $\Delta$, we define its number $j(A)$ of adjacent sites
and its number $\ell(A)$ of loops. In this paper, a loop is still defined by
two adjacent sites capped by another site. Note that this definition is
different from the one used by Bousquet-Mélou \cite{bousquet}, who found
similar results.

With our methods, we have been unable to address the site perimeter, which is
not surprising as the generating function of animals marked with a neighbour
is believed to be non-algebraic on the unbounded lattices \cite{conway}. The
best we could do is to compute bounds on the perimeter, using manipulations
similar to Lemma~\ref{lemma:perimeter}, although we do not give further
details in this paper.

As on the square lattice, we begin by defining generating functions, which are
analogues to the ones of Definitions \ref{definition:square1},
\ref{definition:square2} and \ref{definition:square3}.

\begin{definition}
Let $S\subseteq Q$ be an aligned set of positions. Define the following
generating functions, counting animals with a source included in $S$ (or
inflated heaps with base included in $S$):
\begin{itemize}
\item $\cJ_{[S]}(t)$ and $\cL_{[S]}(t)$, the generating functions of animals
with a source included in $S$, marked respectively with two adjacent sites and
a loop;
\item $\cM_{[S]}(t)$ the generating function of inflated heaps marked with a
piece at a position $q$ such that $q+2\in S$;
\item $\cW_{[S]}(t)$ the generating function of inflated heaps with a minimal
piece at position $q$ and an illegal minimal piece at position $q+2$;
\item $\cI^2_{[S]}(t)$ (resp.\ $\cI^3_{[S]}(t)$) the generating function of
inflated heaps marked with two independent pieces at positions $q$ and $q+2$
(resp.\ $q$ and $q+3$).
\end{itemize}
Also let $\cJ_S(t)$, $\cL_S(t)$, $\cM_S(t)$ and $\cW_S(t)$ be the analogous
generating functions counting animals of source $S$ and heaps of base $S$.
\end{definition}

\begin{theorem}\label{theorem:triang}
The generating functions counting the total number of adjacent sites and loops
of animals of a given area on the triangular lattice are:
\begin{align}\label{js}
\cJ_{[S]}(t)&=\frac{t\cM_{[S]}(t)-\cW_{[S]}(t)}{1+t}\text;\\\label{ls}
\cL_{[S]}(t)&=t\cJ_{[S]}(t)\text;\\\label{jS}
\cJ_S(t)&=\frac{t\cM_S(t)+j(S)\cA_S(t)-\cW_S(t)}{1+t}\text;\\\label{lS}
\cL_S(t)&=t\cJ_S(t)\text.
\end{align}
\end{theorem}

To prove this theorem, we use the results of Section~\ref{subsection:proof},
along with an additional bijection.

\begin{lemma}\label{lemma:psi}
The following identity holds:
\[\cI^2_{[S]}(t)-\cJ_{[S]}(t)=
\frac{t}{1-t}\Bigl(2\cJ_{[S]}(t)+\cI^3_{[S]}(t)\Bigr).\]
\end{lemma}

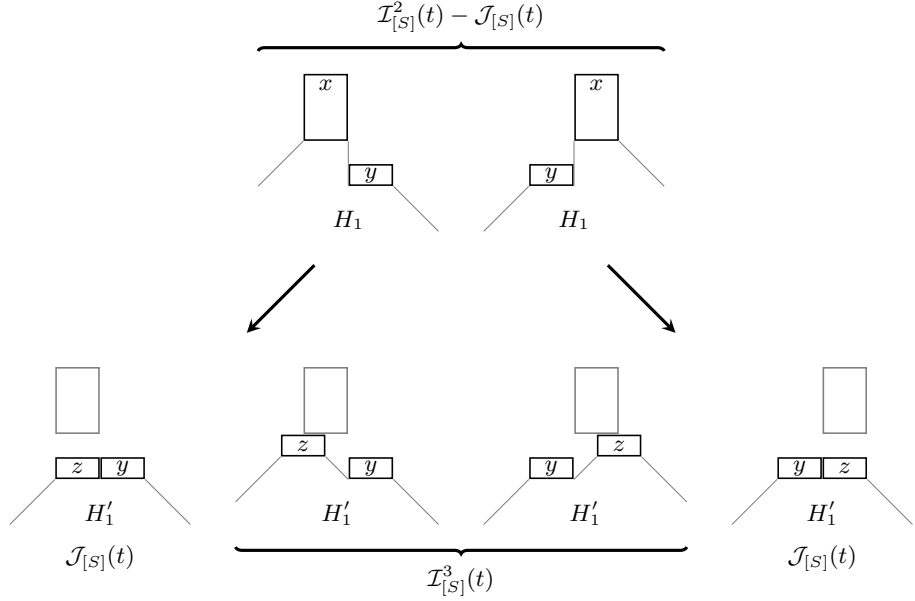
\begin{figure}[ht]
\begin{center}
\begin{tikzpicture}[dominoes]
\begin{scope}
\node at (-1,4) {$x$};
\node at (-1,3) [stack] {};
\node at (1,0) [domino] {$y$};
\draw [help lines] (-2,1.5) -- ++(-2,-2);
\draw [help lines] (2,-.5) -- ++(2,-2);
\draw [help lines] (0,1.5) -- ++(0,-2);
\node at (0,-2) {$H_1$};
%\draw [arrow] (0,-4) -- ++(0,-3);
\end{scope}
\begin{scope}[yshift=-13cm]
\node at (-1,3) [stack, help lines] {};
\node at (-2,1) [domino] {$z$};
\node at (1,0) [domino] {$y$};
\draw [help lines] (-3,.5) -- ++(-2,-2);
\draw [help lines] (2,-.5) -- ++(2,-2);
\draw [help lines] (-1,.5) -- ++(1,-1);
\node at (-.5,-2) {$H_1'$};
\end{scope}
\begin{scope}[xshift=10cm]
\node at (-1,0) [domino] {$y$};
\node at (1,4) {$x$};
\node at (1,3) [stack] {};
\draw [help lines] (-2,-.5) -- ++(-2,-2);
\draw [help lines] (2,1.5) -- ++(2,-2);
\draw [help lines] (0,1.5) -- ++(0,-2);
\node at (0,-2) {$H_1$};
%\draw [arrow] (0,-4) -- ++(0,-3);
\end{scope}
\begin{scope}[xshift=10cm,yshift=-13cm]
\node at (-1,0) [domino] {$y$};
\node at (1,3) [stack, help lines] {};
\node at (2,1) [domino] {$z$};
\draw [help lines] (-2,-.5) -- ++(-2,-2);
\draw [help lines] (3,.5) -- ++(2,-2);
\draw [help lines] (1,.5) -- ++(-1,-1);
\node at (.5,-2) {$H_1'$};
\end{scope}
\begin{scope}[xshift=-11cm,yshift=-13cm]
\node at (-1,0) [domino] {$z$};
\node at (1,0) [domino] {$y$};
\node at (-1,3) [stack, help lines] {};
\draw [help lines] (-2,-.5) -- ++(-2,-2);
\draw [help lines] (2,-.5) -- ++(2,-2);
\node at (0,-2) {$H_1'$};
\end{scope}
\begin{scope}[xshift=21cm,yshift=-13cm]
\node at (-1,0) [domino] {$y$};
\node at (1,0) [domino] {$z$};
\node at (1,3) [stack, help lines] {};
\draw [help lines] (-2,-.5) -- ++(-2,-2);
\draw [help lines] (2,-.5) -- ++(2,-2);
\node at (0,-2) {$H_1'$};
\end{scope}
\draw [brace] (-4,5.5) -- ++(18,0);
\node at (5,7) {$\cI^2_{[S]}(t)-\cJ_{[S]}(t)$};
\draw [brace] (15,-16.5) -- ++(-20,0);
\node at (5,-18) {$\cI^3_{[S]}(t)$};
\node at (-11,-17) {$\cJ_{[S]}(t)$};
\node at (21,-17) {$\cJ_{[S]}(t)$};
\draw [arrow] (-1.5,-4) -- ++(-3,-3);
\draw [arrow] (11.5,-4) -- ++(3,-3);
\end{tikzpicture}
\end{center}
\caption{The bijection $\Psi$: we remove all pieces of the stack of $x$ which
are higher than $y$. This uncovers a piece $z$, either adjacent with $y$ or
with a position at distance~3 and higher than $y$.}\label{figure:psi}
\end{figure}

\begin{proof}
We use a bijection $\Psi$, analogue to $\Phi_1$ (see Lemma~\ref{lemma:phi1})
and illustrated in Figure~\ref{figure:psi}. Let $(H,\{x,y\})$ be a marked
inflated heap, such that $x$ and $y$ are independent, with positions at
distance~2, and not at the same height (say, $x$ is higher). We use the
bijection $\mathcal F_\downarrow$ to form a factorized heap $(H_1\cdot H_2)$.
We then remove from $H_1$ all pieces of the stack of $x$ that are higher than
$y$, thus forming the heap $H_1'$. There are two possiblities:
\begin{enumerate}[label=(\alph*)]
\item $H_1'$ has two maximal pieces, $y$ and $z$, which are adjacent;
\item $H_1'$ has two maximal pieces, $y$ and $z$, with positions at
distance~3.
\end{enumerate}
The inverse bijection is done by putting back the stack of $x$, which can have
an arbitrary number of pieces. In case~(b), as the inflated heap $H_1'$ is
aligned, the pieces $y$ and $z$ cannot be at the same height. Therefore, $z$
must be the higher maximal piece. In case~(a), however, $z$ can be either the
left maximal piece or the right, leading to the factor~2 on the term
$\cJ_{[S]}(t)$.
\end{proof}

\begin{proof}[Proof of Theorem~\ref{theorem:triang}]
First, we derive the identities \eqref{ls} and \eqref{lS} dealing with loops,
using a method identical to the proof of Lemma~\ref{lemma:phi0}. When
dealing with general heaps, there is no $1+t$ factor due to the duplication of
marked pieces.

We now prove the identity \eqref{js}. Let $(H,\{x\})$ be a heap counted by
$\cM_{[S]}(t)$. We use the bijection $\mathcal F_\downarrow$ to pull downwards
the piece $x$, creating a factorized heap. We remark that such factorized
heaps may be built by replacing each piece of an almost strict factorized heap
by an arbitrary stack, leading to the link:
\[\cM_{[S]}(t)=M^*_{[S]}\biggl(\frac{t}{1-t}\biggr)\text.\]
The generating functions $\cI^2_{[S]}(t)$ and $\cI^3_{[S]}(t)$ are also given
in this manner. As for $\cW_{[S]}(t)$, it satisfies:
\[\cW_{[S]}(t)=W_{[S]}\biggl(\frac{t}{1-t}\biggr)
=(1-t)W^*_{[S]}\biggl(\frac{t}{1-t}\biggr).\]
Taking the identities of Lemmas \ref{lemma:phi2} and \ref{lemma:phi3}
together and performing the substitution $t\mapsto\frac{t}{1-t}$, we thus
find:
\[\cI^2_{[S]}(t)+\cW_{[S]}(t)
=\frac{t}{1-t}\Bigl(\cM_{[S]}(t)+\cI^3_{[S]}(t)\Bigr).\]
Using now Lemma~\ref{lemma:psi}, this boils down to \eqref{js}. Performing the
same substitution on the identity of Lemma~\ref{lemma:inclusion} yields:
\[\cW_{[S]}(t)=\sum_{T\subseteq S}\bigl(\cW_T(t)-j(T)\cA_T(t)\bigr).\]
The last identity \eqref{jS} is thus derived using an inclusion-exclusion
argument.
\end{proof}

\section{Asymptotic results}\label{section:asymptotic}

\subsection{Animals according to area}

Here, we derive asymptotic estimates from the results of
Section~\ref{section:dominoes}. First, consider the polynomials $F_m(t)$ and
$\widehat F_m(t)$, defined in Definition~\ref{definition:fm}. Let $\rho_m$ and
$\sigma_m$ be their respective smallest root.

\begin{lemma}\label{lemma:roots}
For all $m\geq0$, the polynomials $F_m(t)$ and $\widehat F_m(t)$ have only
real, simple roots. Their smallest roots $\rho_m$ and $\sigma_m$ verify:
\begin{align*}
\frac1{\rho_m}&=4\cos^2\frac\pi{m+2}\text;\\
\frac1{\sigma_m}&=4\cos^2\frac\pi{2m}\text.
\end{align*}
\end{lemma}

\begin{proof}
We check by induction on $m$ that the degrees of $F_m(t)$ and $\widehat
F_m(t)$ are \smash{$\bigl\lceil\frac m2\bigr\rceil$} and
\smash{$\bigl\lfloor\frac m2\bigr\rfloor$}, respectively. We also check by
induction the following identities:
\begin{align*}
F_m\biggl(\frac1{4\cos^2\theta}\biggr)&=
\frac{\sin[(m+2)\theta]}{(2\cos\theta)^{m+1}\sin\theta}\text;\\
\widehat F_m\biggl(\frac1{4\cos^2\theta}\biggr)&=
\frac{2\cos(m\theta)}{(2\cos\theta)^{m}}\text.
\end{align*}
By choosing appropriate values of $\theta$ in the interval $[0,\pi/2)$, these
identities account for all the roots of the polynomials. Thus, we prove that
the roots are real and simple, and we derive the value of the smallest root.
\end{proof}

Now, let $\Gamma$ be a square lattice and $\Delta$ be its associated
triangular lattice; let $S$ be a one-line source. We denote by $a(n)$ and
$\bar a(n)$ the number of animals of area~$n$ of source $S$ on the lattices
$\Gamma$ and $\Delta$ respectively. The result below is simply obtained by
performing singularity analysis \cite{flajolet} on the formul{\ae} of
Section~\ref{section:dominoes}.

\begin{proposition}\label{prop:asympt}
The general form of the asymptotic behaviour of $a(n)$ and $\bar a(n)$ is:
\begin{align*}
a(n)&\sim\lambda\mu^nn^\nu\text;&
\bar a(n)&\sim\bar\lambda\bar\mu^nn^\nu\text,
\end{align*}
where the constants $\bar\mu$ and $\nu$ are:
\begin{itemize}
\item in the full lattice, $\bar\mu=4$ and $\nu=-1/2$;
\item in the half lattice, $\bar\mu=4$ and $\nu=-3/2$;
\item in the cylindrical lattice of width~$m$, $\bar\mu=1/\sigma_m$ and $\nu=0$;
\item in the rectangular lattice of width~$m$, $\bar\mu=1/\rho_m$ and $\nu=0$.
\end{itemize}
Moreover, in each case, $\mu$ is equal to $\bar\mu-1$ and $\lambda$ and
$\bar\lambda$ depend on the source $S$.
\end{proposition}

Notably, changing the source $S$ only changes the behaviour of $a(n)$ and
$\bar a(n)$ by a multiplicative constant.

\subsection{Average number of adjacent sites and loops and average perimeter}

Let now $j(n)$ be the average number of adjacent sites in the animals of
source~$S$ and area~$n$ in the lattice $\Gamma$. Let $\ell(n)$, $p_e(n)$,
$p_i(n)$ be their average number of loops, external perimeter, and internal
perimeter; let $\bar\jmath(n)$ and $\bar\ell(n)$ be the analogous quantities
in the lattice $\Delta$.

\begin{corollary}
Assume that $\Gamma$ is either the full lattice, the half lattice, or a
cylindrical bounded lattice. As $n$ tends to infinity, we have the following
estimates:
\begin{align*}
j(n)&\sim\frac n{\mu+1};&\bar\jmath(n)&\sim\frac n{\bar\mu+1}\text;\\
\ell(n)&\sim\frac n{\mu^2};&\bar\ell(n)&\sim\frac n{\bar\mu(\bar\mu+1)}\text;\\
p_i(n)&\sim p_e(n)\sim\frac\mu{\mu+1}n\text.
\end{align*}
\end{corollary}

In the unbounded lattices, the growth constants are $\bar\mu=4$ and $\mu=3$.
Thus, these estimates become:
\begin{align*}
j(n)&\sim\frac n4;&\ell(n)&\sim\frac n9;&p(n)&\sim\frac{3n}4;&
\bar\jmath(n)&\sim\frac n5;&\bar\ell(n)&\sim\frac n{20}.
\end{align*}

\begin{proof}
Let us begin with the number of adjacent pieces in the square lattice. This
number is given by the identity \eqref{JS}:
\[J_S(t)=\frac{tM_S(t)+j(S)A_S(t)-W_S(t)}{1+t},\]
where the generating functions are defined in Definitions
\ref{definition:square1} and \ref{definition:square2}. We examine the
coefficients of these generating functions.
\begin{itemize}
\item The $n$th coefficient of $J_S(t)$ is $j(n)a(n)$.
\item As, in the full, half, and cylindrical models, the position $q+2$ is
always in $Q$ as soon as $q$ is, the generating function $M_S(t)$ simply
counts animals marked with any site; its $n$th coefficient is $na(n)$.  \item
As a corollary to Proposition~\ref{prop:asympt}, the $n$th coefficient of both
$A_S(t)$ and $W_S(t)$ is \smash{$\mathrm O\bigl(a(n)\bigr)$}.
\end{itemize}
From this, it follows that the dominant term in the right-hand side is
that of $M_S(t)$. We perform singularity analysis, letting $t$ tend to
the singularity $1/\mu$. We obtain, as $n$ tends to infinity:
\[j(n)a(n)\sim\frac{1/\mu}{1+1/\mu}na(n).\]
The result follows; the other estimates are obtained with a similar analysis
on the equations of Theorems \ref{theorem:square} and \ref{theorem:triang}.
\end{proof}

\section{Examples}\label{section:examples}

\subsection{Single-source animals on the full lattice}

We start with the simplest case, that of single-source animals on the full
lattices.

\begin{corollary}
The generating functions counting the total number of adjacent sites, number
of loops and site perimeter of the single-source directed animals on the full
square lattice are respectively given by:
\begin{align*}
J(t)&=\frac1{2t(1+t)}
\biggl(1-\frac{1-4t+t^2+4t^3}{\sqrt{1+t}(1-3t)^{3/2}}\biggr)\text;\\
L(t)&=\loops\text;\\
P(t)&=\conway\text.
\end{align*}
\end{corollary}

The value of $P(t)$ was conjectured by Conway \cite{conway}, and the value of
$L(t)$ was proved by Bousquet-Mélou using a gas model method \cite{bousquet}.

\begin{proof}
We use Theorem~\ref{theorem:square} to derive the generating functions. First,
we use \eqref{JS} to compute $J(t)\equiv J_{\{0\}}(t)$, which gives:
\[J_{\{0\}}(t)=
\frac{tM_{\{0\}}(t)+j\bigl(\{0\}\bigr)A_{\{0\}}(t)-W_{\{0\}}(t)}{1+t}.\]
The generating function $M_{\{0\}}(t)$ is simply equal to $tA'(t)$,
$j(\{0\})$ is zero, and $W_{\{0\}}(t)$ is equal to $A_{\{0,2\}}(t)$, in
turn equal to $D(t)A(t)$ using Proposition~\ref{proposition:compact}. This
yields the announced formula; the other two generating functions follow from
equations \eqref{LS} and \eqref{PeS}.
\end{proof}

Similarly, Theorem~\ref{theorem:triang} instantiates on single-source animals
on the triangular full lattice. We omit the proof, which is identical to the
square lattice case.

\begin{corollary}
The generating functions counting the total number of adjacent sites and number of loops of single-source animals on the full triangular lattice are:
\begin{align*}
\cJ(t)
&=\frac1{2t(1+t)}\biggl(1-t-\frac{1-7t+12t^2-2t^3}{(1-4t)^{3/2}}\biggr)\text;\\
\cL(t)
&=\frac1{2(1+t)}\biggl(1-t-\frac{1-7t+12t^2-2t^3}{(1-4t)^{3/2}}\biggr)\text;\\
\end{align*}
\end{corollary}

This time, the value of $\mathcal L(t)$ is different from the one found by
Bousquet-Mélou \cite{bousquet}, who used a different definition of loops.

\subsection{Compact-source animals on the full lattice}

As an illustration of how to deal with non-single source animals, we consider
animals with any compact source (see Definition~\ref{definition:compact}).
Recall that the number of such animals of area~$n$ is $3^{n-1}$ on the square
full lattice, so that the generating function is:
\[A_c(t)=\frac{t}{1-3t}.\]
We only give the result for the number of adjacent sites of animals on the
square full lattice, but other configurations can be handled similarly.

\begin{corollary}
The generating function counting the total number of adjacent sites of the
compact-source directed animals on the full square lattice is:
\begin{equation*}
J_c(t)=\frac12\Biggl(
\frac{1-2t}{\sqrt{1+t}(1-3t)^{3/2}}-\frac{1-3t-2t^2}{(1+t)(1-3t)^2}\Biggr).\\
\end{equation*}
\end{corollary}

\begin{proof}
Let $C$ be a compact source with $k$ sites. The generating function $A_C(t)$
is, according to Proposition~\ref{proposition:compact}:
\[A_C(t)=D(t)^{k-1}A(t).\]
Moreover, $M_C(t)$ is simply equal to $tA_C'(t)$, $j(C)$ is $k-1$, and
$W_C(t)$ counts animals with a compact source with $k+1$ sites, and is thus
equal to $D(t)A_C(t)$. Therefore:
\[J_C(t)=\frac{t^2A_C'(t)+(k-1)A_C(t)-D(t)A_C(t)}{1+t}.\]
We sum this identity for all $k\geq0$:
\[J_c(t)=\frac1{1+t}\Biggl(t^2A_c'(t)
+\frac{A(t)D(t)}{\bigl(1-D(t)\bigr)^2}
-D(t)A_c(t)\Biggr).\]
This completes the proof.
\end{proof}

\subsection{Half-animals on the square rectangular lattices}

Finally, we present our results on the external and internal site perimeter of
half-animals (that is, animals of source~$\{0\}$) on the square rectangular
lattices. The former was the object of a conjecture by Le~Borgne
\cite{leborgne}; from our formula, one can prove this conjecture.

We denote by $D_m(t)$ the generating function of half-animals in square
rectangular lattice of width~$m$.

\begin{corollary}
The generating functions giving the total external and internal site perimeter
of half-animals on the square rectangular lattice of width~$m$ are:
\begin{align*}
P^e_m(t)&=D_m(t)+\frac{t}{1+t}D_m'(t)+\frac1{1+t}D_m(t)^2\text;\\
P^i_m(t)&=\frac{t}{1+t}D_m(t)+\frac{t}{1+t}D_m'(t)
-\frac1{1+t}D_m(t)\Bigl(D_m(t)-D_{m-2}(t)\Bigr)\text,
\end{align*}
where the generating function $D_m(t)$ is derived from \eqref{dm}.
\end{corollary}

By letting $m$ tend to infinity, we obtain the generating functions of animals
on the half lattice:
\begin{align*}
P^e(t)&=D(t)+\frac{t}{1+t}D'(t)+\frac1{1+t}D(t)^2;\\
P^i(t)&=\frac{t}{1+t}D(t)+\frac{t}{1+t}D'(t).
\end{align*}

\begin{proof}
Let $Q$ be the rectangular model $\{0,\dotsc,m-1\}$ of width~$m$. By symmetry,
instead of considering animals with a source at position~$0$, we consider them
to have a source at position~$m-1$. This does not change the site perimeter of
the animals.

The generating functions $P^e_m(t)\equiv P^e_{\{m-1\}}(t)$ and $P^i_m(t)\equiv
P^i_{\{m-1\}}(t)$ are given by \eqref{PeS} and \eqref{PiS}, which in turn
require us to compute the generating function $J_{\{m-1\}}(t)$. The number
$j(\{m-1\})$ is again zero; moreover, as the position $m+1$ is not in $Q$,
$W_{\{m-1\}}(t)$ is also zero. Thus, all we need to compute are the generating
functions $M_{\{m-1\}}(t)$ and $E_{\{m-1\}}(t)$.

Let $D_m^{(q)}(t)$ be the generating function of half-animals marked with a
site at position $q$. Using Definition~\ref{definition:square1}, we find:
\begin{align*}
M_{\{m-1\}}(t)&=tD_m'(t)-D_m^{(m-1)}(t)-D_m^{(m-2)}(t);\\
E_{\{m-1\}}(t)&=D_m^{(m-1)}(t)+D_m^{(0)}(t).
\end{align*}

Finally, we derive the generating functions $D_m^{(q)}(t)$ using
Lemma~\ref{lemma:marked}; in the notations of this lemma, $D_m^{(q)}(t)$ is
equal to $H_{[\{m-1\}]}^{(q)}(t)$ as a marked heap cannot be empty. We must
compute the following generating functions:
\begin{itemize}
\item $H_{\{m-1\}}(t)$ and $H_{\{0\}}(t)$ are both equal to $D_m(t)$ by
symmetry;
\item the only heap of base included in $\{m-1\}$ avoiding $m-2$ is the empty
heap, so that $V^{m-2}_{[\{m-1\}]}(t)=1$;
\item as a strict pyramid of base~$m-1$ is either a single piece or a piece
topped by a pyramid of base~$m-2$, we have
\smash{$H_{\{m-2\}}(t)=\frac{D_m(t)}t-1$};
\item as a pyramid of base $m-1$ avoiding $0$ lives in the model
$\{2,\dotsc,m-1\}$, which has $m-2$ positions, we have
$V^0_{[\{m-1\}]}(t)=1+D_{m-2}(t)$.
\end{itemize}
From this, we find:
\begin{align*}
D_m^{(m-1)}(t)&=\frac1{1+t}\Bigl(1+D_m(t)\Bigr)D_m(t);\\
D_m^{(m-2)}(t)&=\frac1{1+t}D_m(t)\biggl(\frac{D_m(t)}t-1\biggr);\\
D_m^{(0)}(t)&=\frac1{1+t}\Bigl(D_m(t)-D_{m-2}(t)\Bigr)D_m(t).
\end{align*}
Injecting these values into \eqref{JS}, \eqref{PeS} and \eqref{PiS}, we get
the announced results.
\end{proof}

\section*{Acknowledgments}
I would like to thank Mireille Bousquet-Mélou for her precious help in the
writing of this paper.

\bibliographystyle{hplain}
\bibliography{animals}{}

\end{document}